\documentclass[a4paper,12pt,twoside]{article}

\usepackage[dvipsnames]{xcolor}

\usepackage{float} 

\usepackage{hyperref} %Internal and external links
\hypersetup{breaklinks, colorlinks,
    linkcolor = {Blue},
    citecolor = {Blue},
    urlcolor  = {Blue}
}

\usepackage{amsmath,amssymb}

\usepackage[text={6.5in,9.5in},centering]{geometry}

\usepackage{theorem}

\usepackage[inline,shortlabels]{enumitem} 
\renewlist{enumerate}{enumerate}{1}
\setlist[enumerate,1]{label={\upshape({\alph*})},ref={\alph*},align=left,
labelsep=0.5ex,leftmargin=*}

\DeclareSymbolFont{bbold}{U}{bbold}{m}{n}
\DeclareMathSymbol{\bbone}{\mathord}{bbold}{49}    % mathbb 1

\newtheorem{theorem}{Theorem}[section]
\newtheorem{corollary}{Corollary}[section]
\newtheorem{lemma}{Lemma}[section]

{\theorembodyfont{\upshape}
\newtheorem{remark}{Remark}[section]
\newtheorem{example}{Example}[section]
}

\newcommand{\Proof}{\textbf{Proof. }}

\newcommand{\CC}{\mathbb{C}}
\newcommand{\RR}{\mathbb{R}}
\newcommand{\NN}{\mathbb{N}}
\newcommand{\ZZ}{\mathbb{Z}}

\newcommand{\calz}{\mathcal{Z}}

\newcommand{\ee}{{\mathrm e}}

\newcommand{\card}[1]{\lvert#1\rvert}       % cardinality of a set

\newcommand{\citet}{\cite}

\newcommand{\Coeff}[2]{\mathrm{Coeff}\Big(#1;#2\Big)}

\newcommand{\coeffw}{\mathrm{Coeff}}

\renewcommand{\geq}{\geqslant}

\newcommand{\leer}[1]{}

\renewcommand{\leq}{\leqslant}

\newcommand{\Per}{\mathrm{Per}}

\newcommand{\qed}{$\square$}                % end of proof
 
\newcommand{\set}[1]{\underline{#1}}

\newcommand{\setn}[1]{\underline{#1}_0} 
 
\newcommand{\binomial}[2]{\genfrac{(}{)}{0pt}{}{#1}{#2}}

\newcommand{\newatop}[2]{\genfrac{}{}{0pt}{}{\scriptstyle #1}{\scriptstyle #2}}

\scrollmode
\begin{document}
%%%%%%%%%%%%%%%%%%%%%%%%%%%%%%%%%%%%%%%%%%%%%%%%%%%%%%%%%%%%%%%%%%%%%%
%%%%%%%%%%%%%%%%%%%%%%%%%%%%%%%%%%%%%%%%%%%%%%%%%%%%%%%%%%%%%%%%%%%%%%
\title{
New permanent approximation inequalities via identities
}
\author{Bero Roos\footnote{
Postal address: FB IV -- Department of Mathematics,
University of Trier,
54286 Trier,
Germany. 
E-mail address: \texttt{bero.roos@uni-trier.de}}\\
University of Trier
}
\date{ \today }
\maketitle
%%%%%%%%%%%%%%%%%%%%%%%%%%%%%%%%%%%%%%%%%%%%%%%%%%%%%%%%%%%%%%%%%%%%%%
%%%%%%%%%%%%%%%%%%%%%%%%%%%%%%%%%%%%%%%%%%%%%%%%%%%%%%%%%%%%%%%%%%%%%%
\begin{abstract}
The aim of this paper is to present new upper bounds for the 
distance between a properly normalized permanent of a rectangular 
complex matrix and the product of the arithmetic means of the 
entries of its columns. It turns out that the bounds 
improve on those from earlier work. 
Our proofs are based on some new identities for 
the above-mentioned  differences and also for related expressions
for matrices over a rational associative commutative unital 
algebra. Some of our identities are generalizations of results 
in Dougall (Proc. Edinburgh Math. Soc.,
24:61--77, 1905). Second order results are also included. 
\medskip\\
%%%%%%%%%%%%%%%%%%%%%%%%%%%%%%%%%%%%%%%%%%%%%%%%%%%%%%%%%%%%%%%%%%%%%%
{\small 2010 Mathematics Subject Classification:
Primary
15A15;    % Determinants, permanents, other special matrix functions
secondary
15A45;    % Miscellaneous inequalities involving matrices
05A99.    % Combinatorics: None of the above, but in this section
}

\end{abstract}
%%%%%%%%%%%%%%%%%%%%%%%%%%%%%%%%%%%%%%%%%%%%%%%%%%%%%%%%%%%%%%%%%%%%%%
%%%%%%%%%%%%%%%%%%%%%%%%%%%%%%%%%%%%%%%%%%%%%%%%%%%%%%%%%%%%%%%%%%%%%%
\section{Introduction} 
%%%%%%%%%%%%%%%%%%%%%%%%%%%%%%%%%%%%%%%%%%%%%%%%%%%%%%%%%%%%%%%%%%%%%%
%%%%%%%%%%%%%%%%%%%%%%%%%%%%%%%%%%%%%%%%%%%%%%%%%%%%%%%%%%%%%%%%%%%%%%
It is well-known that computing the permanent of an 
$n\times n$ matrix can be a difficult task, if $n$ is a large 
natural number, see Valiant \citet{MR526203} and 
Minc \citet[Chapter 7]{MR504978}. 
There are a couple of known explicit formulae, 
the most efficient of which
seem to be due to  Ryser \citet[Theorem 4.1, page 26]{MR0150048} or 
Glynn \citet[Theorem 2.1]{MR2673027} and require at least $O(2^nn)$ 
arithmetic operations. 
Matrices with a special structure 
can sometimes be treated differently, 
e.g.\ see Minc \citet[Section 3.4 or Lemma 1 on page 113]{MR504978},
Bax and Franklin \citet{MR1470132},  Schwartz \citet{MR2489398},
Bj\"orklund et al.\ \citet{MR2682986}  and the references therein. 
On the other hand, there are approximation algorithms, 
e.g.\ see Jerrum et al.\ \citet{MR2147852}, 
Barvinok \citet{MR3464208} and the references 
given there. 

There are many upper and lower bounds for permanents, see e.g.\ 
Minc \citet[Chapters 4--6]{MR504978}. 
But the literature seems to contain only a few 
explicit approximation inequalities. 
See Bhatia \citet{MR761074}, Elsner \citet{MR961564}, 
Bhatia and Elsner \citet{MR1054132}, 
and Friedland \citet{MR1054131}, for some upper bounds of the distance between 
two permanents of quadratic complex matrices; however,
the inequalities given there are not easily comparable
with those discussed below. 
Another approach is to approximate a permanent by more special 
expressions, which depend on the matrix under consideration. 
Here, for convenience, we consider properly normalized permanents, 
i.e.\ we divide by the number of summands. 

In this paper, we consider the approximation of a normalized 
permanent of a rectangular complex matrix. 
If the rows of the matrix are approximately equal, a 
good approximant should be the product of the arithmetic means of the 
entries of the columns of this matrix.
Some upper bounds for the approximation error can be found in 
Bobkov \citet{MR2137450} and Roos \citet{MR3314090}.

To be more precise, we need the following notation. 
Let $N\in\NN$, $n\in\set{N}:=\{1,\dots,N\}$ and 
$Z=(z_{j,r})\in\CC^{N\times n}$ be an $N\times n$ matrix with
complex entries. We set 
$\widetilde{z}_r=\frac{1}{N}\sum_{j=1}^N z_{j,r}$, $(r\in\set{n})$ 
and assume that $|z_{j,r}|\leq 1$, $(j\in \set{N},\,r\in\set{n})$. 
However, it is noteworthy that some of the results of 
Section \ref{s749687} below do not require the latter boundedness 
assumption. 

For arbitrary sets $A$ and $B$, let $A^B$, resp.\ 
$A_{\neq}^B$, be the set of all maps, resp.\ injective maps,  
$f:\,B\longrightarrow A$. 
For $f\in A^B$ and $b\in B$, we write
$f(b)=f_b$. Let $\set{N}^n=\set{N}^{\set{n}}$ and 
$\set{N}_{\neq}^n=
\set{N}_{\neq}^{\set{n}}
=\{(j_1,\dots,j_n)\in\set{N}^n\,|\,j_r\neq j_s
\mbox{ for all }r,s\in\set{n}\mbox{ with }r\neq s\}$. 
In particular, $\set{N}_{\neq}^{N}$ is the set of all permutations 
on the set $\set{N}$. The permanent of $Z$ can now be defined by
\begin{align*}
\Per(Z)=\sum_{j\in \set{N}_{\neq}^n}\prod_{r=1}^n z_{j_r,r}.
\end{align*}
As indicated above,
\begin{align*} 
\frac{(N-n)!}{N!}\Per(Z)\approx\prod_{r=1}^n\widetilde{z}_r,
\end{align*}
when 
\begin{align}\label{e621850987}
z_{1,r}\approx\dots\approx z_{N,r}\mbox{ for all }r\in\set{n}. 
\end{align}
We note that, if $Z$ has identical columns, i.e.\ 
$z_{j,1}=\dots= z_{j,n}$ for all 
$j\in\set{N}$, then we have 
$\prod_{r=1}^n\widetilde{z}_r=\widetilde{z}_1^n$,
whereas
\begin{align*}
\frac{(N-n)!}{N!}\Per(Z)=\frac{1}{\binomial{N}{n}}
\sum_{\newatop{J\subseteq \set{N}}{\card{J}=n}}\prod_{j\in J}z_{j,1}
\end{align*}
is the normalized elementary symmetric polynomial of degree $n$
in the variables $z_{1,1},\dots, z_{N,1}$.
Here, for a finite set $J$,  let $\card{J}$ be the number of
its elements.

Let us give a review of some approximation inequalities from the 
literature. Bobkov
\citet[Theorem 2.1]{MR2137450} showed by a somewhat 
complicated induction that 
\begin{align}\label{e6285}
\Big|\frac{(N-n)!}{N!}\Per(Z)-\prod_{r=1}^n\widetilde{z}_r\Big|
\leq C_0\frac{n}{N} \qquad \mbox{with } C_0=16
\end{align}
and used this inequality to study 
an approximate de Finetti representation for probability measures, 
on product measurable spaces, which are symmetric under permutations 
of coordinates. The upper bound in \eqref{e6285} is small if $n$ is 
small in comparison with $N$.  But since it is independent of $Z$,
it is not good in the case \eqref{e621850987}.

A bound depending on $Z$ was given in Roos \citet{MR3314090}.
From the more general Theorem 2.13 given there, it follows that 
\begin{align}\label{e4176438}
\Big|\frac{(N-n)!}{N!}\Per(Z)-\prod_{r=1}^n\widetilde{z}_r\Big|
\leq 3.57\gamma, 
\end{align}
where 
\begin{gather*}
\gamma=\gamma(1), \quad
\gamma(x)=\frac{n\alpha}{N}\min\Big\{xn,\frac{1}{1-\beta}\Big\},
\qquad (x\in[0,\infty)),\nonumber\\
\alpha=\frac{1}{nN}\sum_{j=1}^N\sum_{r=1}^n|a_{j,r}|^2, \qquad 
a_{j,r}=z_{j,r}-\widetilde{z}_r,\quad 
(j\in\set{N},\,r\in\set{n}),\qquad 
\beta=\frac{1}{n}\sum_{r=1}^n|\widetilde{z}_r|^2.  
\end{gather*}
In Remark 2.9 of that paper, it was also shown that 
$\gamma\leq\frac{n}{N}$. Consequently in \eqref{e6285},
$C_0$ can be replaced with $3.57$. However, inequality 
\eqref{e4176438} is preferable to \eqref{e6285} with 
any constant $C_0$, since $\gamma$ can be much smaller than 
$\frac{n}{N}$. In fact, the right-hand side in \eqref{e4176438}
is small in the case \eqref{e621850987}.

The proof of \eqref{e4176438} does not require an induction 
argument but instead is based on the representation 
(see \cite[Theorem 2.8]{MR3314090})
\begin{align*}
\frac{(N-n)!}{N!}\Per(Z)
=H_n(Z),
\end{align*}
where $H_\ell(Z)=\sum_{m=0}^\ell G_m(Z)$ for $\ell\in\set{n}$, 
\begin{align*}
G_m(Z)
&=\frac{(N-m)!}{(n-m)!N!}
\Coeff{x_1\cdots x_n}{
\Bigl(\sum_{r=1}^n\widetilde{z}_rx_r\Bigr)^{n-m}
\prod_{j=1}^N\Bigl(1+\sum_{r=1}^na_{j,r}x_r\Bigr)},
\end{align*}
for $m\in\setn{n}=\{0,\dots,n\}$, 
and $\coeffw$ denotes the coefficient of $x_1\cdots x_n$ 
in the formal power series expansion of the expression given above. 
In particular, 
$H_1(Z)=\prod_{r=1}^n\widetilde{z}_r$, and, if $n\geq2$,  
\begin{align*}
H_2(Z)
=\prod_{r=1}^n\widetilde{z}_r-\frac{1}{N(N-1)}
\sum_{\newatop{R\subseteq\set{n}}{\card{R}=2}}
\Bigl(\sum_{j=1}^N\prod_{r\in R}a_{j,r}\Bigr)
\prod_{r\in\set{n}\setminus R}\widetilde{z}_r.
\end{align*}
It turned out that $\frac{(N-n)!}{N!}\Per(Z)$ can be approximated 
by $H_\ell(Z)$, $(\ell\in\set{n})$, which we call the $\ell$th 
order approximant. In fact, the following estimate 
shows that the accuracy is increasing in $\ell$: if $\gamma<1$, then
\begin{align*}
\Big|\frac{(N-n)!}{N!}\Per(Z)-H_\ell(Z)\Big|
\leq (\ell+1)^{1/4}\widetilde{C}_{\ell+1}
\frac{\gamma^{(\ell+1)/2}}{(1-\gamma)^{3/4}}, 
\end{align*}
where 
$\widetilde{C}_\ell=(\frac{\ee^\ell\,\ell!}{\ell^{\ell+1/2}})^{1/2}$.
We note that Corollary 2.12 in \cite{MR3314090}
gives in the case $\ell\in\set{2}$ and $\gamma<1$ 
the sometimes sharper bounds for the first and second order 
approximations:
\begin{gather} %  number of formula is necessary !
\Big|\frac{(N-n)!}{N!}\Per(Z)-\prod_{r=1}^n\widetilde{z}_r\Big|
\leq \gamma(1/2)+\frac{2.12\,\gamma^{3/2}}{(1-\gamma)^{3/4}},
\label{e63286508}\\
\Big|\frac{(N-n)!}{N!}\Per(Z)-H_2(Z)\Big|
\leq \sqrt{3}\sum_{j=1}^N
\Big(\frac{1}{N^2}\sum_{r=1}^n|a_{j,r}|^2\min\Big\{
\frac{n}{3},\,\frac{1}{1-\beta}\Big\}\Big)^{3/2}
+\frac{2.27\gamma^2}{(1-\gamma)^{3/4}}, \label{e4176439}
\end{gather}
where, for \eqref{e4176439}, we assume that $n\geq 2$.
Hence, if $\gamma$ is small, then 
$|\frac{(N-n)!}{N!}\Per(Z)-\prod_{r=1}^n\widetilde{z}_r|$ 
is bounded by $C_1\gamma(\frac{1}{2})$ with $C_1\approx 1$.

The results of the present paper imply that, in \eqref{e4176438}
or~\eqref{e63286508},
not only the constants but also the form of the right-hand side 
can substantially be improved, see Theorems~\ref{t74739} 
and~\ref{t4216459} below. 
In particular, Theorem \ref{t74739} implies that, 
if $2\leq n\leq N$ and 
\begin{gather*}
\vartheta=\frac{1}{N(N-1)\sqrt{n(n-1)}}
\Big(\sum_{(r,s)\in \set{n}_{\neq}^2}
\Big(\sum_{(u,v)\in\set{N}_{\neq}^2}
|z_{u,r}-z_{v,r}||z_{u,s}-z_{v,s}|\Big)^2\Big)^{1/2},
\end{gather*}
then 
\begin{align}\label{e62865098}
\Big|\frac{(N-n)!}{N!}\Per(Z)-\prod_{r=1}^n\widetilde{z}_r\Big|
&\leq \frac{n-1}{2N}\vartheta
\frac{1-\beta^{n/4}}{1-\sqrt{\beta}},
\end{align}
see \eqref{e739065}. 
Here, the right-hand side of \eqref{e62865098} can be further 
estimated by $(1+\sqrt{\beta})\gamma(\frac{1}{2})
\leq 2\gamma(\frac{1}{2})$, see Remark~\ref{r732860} below. 
However, \eqref{e62865098} can be much better than these alternative 
bounds, see Parts \eqref{ex782566.b} and \eqref{ex782566.c} of 
Example \ref{ex782566} on derangement and m\'{e}nage numbers. 
Indeed, we obtain bounds of 
the order $O(\frac{1}{n})$ and $O(\frac{1}{\sqrt{n}})$ as
$n\to\infty$, whereas the upper bounds in \eqref{e4176438} and 
\eqref{e63286508} cannot be small, since they
contain one of the terms $\gamma$ or $\gamma(\frac{1}{2})$.
The present paper also contains an improvement of 
\eqref{e4176439}, which however is more complicated, 
see Theorem \ref{th487698}.

Let us comment on the method used in this paper. 
Our approach consists of two steps. 
First, we develop some identities for the difference of  
$\Per(Z)$ and its approximant. 
After that, these identities together with the properties of the 
norm and further auxiliary inequalities for permanents 
(see Lemma \ref{l378687}) are applied.
We do not use the methods of \cite{MR2137450} or \cite{MR3314090}. 

Our identities are not only valid for complex matrices, 
but also for matrices over a rational associative commutative 
unital algebra. In the theory of permanents one often
considers matrices over a commutative ring (see Minc
\citet[page 1]{MR504978}), but this is not sufficient here, 
since we need to be able to multiply with rational numbers. 
Some of our identities are generalizations of old identities of 
Dougall \citet{Dougall1905}, who considered, among other things, 
the difference $\prod_{j=1}^Nz_j-\widetilde{z}^N$, where
$z_1,\dots,z_N\in\CC$ and 
$\widetilde{z}=\frac{1}{N}\sum_{j=1}^Nz_{j}$.  
In fact, our first result is Theorem \ref{t6296},
which is a generalization of formula~(3) in 
\cite[page~65]{Dougall1905} concerning elementary 
symmetric polynomials, see Corollary \ref{c76276597} below.
The latter result was a starting point for several other results in 
\cite{Dougall1905}. Similarly, our generalization implies the 
identity (see \eqref{e8366548} below)
\begin{align}\begin{split}
\Per(Z)
-\frac{N!\,}{(N-n)!}\prod_{r=1}^n\widetilde{z}_{r}
&=-\sum_{k=2}^n\frac{1}{2Nk\binomial{n}{k}}
\sum_{\newatop{R\subseteq\set{n}}{\card{R}=k}}
\sum_{(r,s)\in R_{\neq}^2}\sum_{j\in\set{N}_{\neq}^{n}}
(z_{j_r,r}-z_{j_s,r})(z_{j_r,s}-z_{j_s,s})\\
&\qquad{}\times
\Big(\prod_{\ell\in R\setminus\{r,s\}}z_{j_\ell,\ell}\Big)
\prod_{\ell\in\set{n}\setminus R}\widetilde{z}_{\ell},
\end{split}\label{e628507}
\end{align}
which, in turn, is a generalization of another identity in
\cite[page 77]{Dougall1905}, see Corollary \ref{c526745} below. 

We note that, in \eqref{e628507}, it is important to have 
the product $(z_{j_r,r}-z_{j_s,r})(z_{j_r,s}-z_{j_s,s})$ of 
two differences of certain entries of $Z$. 
As a rule, an accurate approximation 
of $\Per(Z)$ should be reflected in a high number of such differences
in the corresponding identity. Indeed, 
Theorem~\ref{t384696} contains an identity for the difference
of $\Per(Z)-\frac{N!}{(N-n)!}H_2(Z)$, where
the right-hand side consists of two expressions containing the 
product of three, resp.\ four, such differences.

The paper is structured as follows. Section \ref{s326765}
is devoted to the notation, which is needed to simplify the 
presentation. In Section~3, we derive  
some new identities for permanents
and related expressions, some of which will 
be used in Section \ref{s749687} to give refined upper bounds of
$|\frac{(N-n)!}{N!}\Per(Z)-H_\ell(Z)|$ for $\ell\in\set{2}$. 
%%%%%%%%%%%%%%%%%%%%%%%%%%%%%%%%%%%%%%%%%%%%%%%%%%%%%%%%%%%%%%%%%%%%%%
\section{Notation}\label{s326765}
%%%%%%%%%%%%%%%%%%%%%%%%%%%%%%%%%%%%%%%%%%%%%%%%%%%%%%%%%%%%%%%%%%%%%%
From now on, unless stated otherwise, our notation is as follows. 
Let $\calz$ be a rational associative commutative unital algebra,
$N\in\NN$, $n\in\set{N}=\{1,\dots,N\}$, 
$Z=(z_{j,r})\in\calz^{N\times n}$, 
\begin{align*}
\widetilde{z}_r=\frac{1}{N}\sum_{j=1}^Nz_{j,r} \mbox{ for } 
r\in\set{n} \quad \mbox{ and }\quad 
y_{j,k,r}=z_{j,r}-z_{k,r} \mbox{ for } j,k\in\set{N},\, r\in\set{n}.
\end{align*}
 Let $p_{j,R}=\prod_{r\in R}z_{j_r,r}$ 
for $R\subseteq\set{n}$ and $j\in\set{N}_{\neq}^S$, 
whenever $S\subseteq \set{N}$ with $R\subseteq S$. 
For $R\subseteq\set{n}$, set 
\begin{align*}
\overline{p}_R=\sum_{j\in\set{N}_{\neq}^{n}} p_{j,R} 
\quad\mbox{ and }\quad 
\widetilde{p}_{R}=\prod_{r\in R}\widetilde{z}_r.
\end{align*}
In particular, we have 
\begin{gather*}
\overline{p}_{\emptyset}=\frac{N!}{(N-n)!},\qquad
\overline{p}_{\{r\}}=\frac{N!}{(N-n)!}\widetilde{z}_r
\quad \mbox{ for }r\in\set{n},\qquad
\overline{p}_{\set{n}}=\Per(Z),\\
\widetilde{p}_{\emptyset}=1,\qquad 
\widetilde{p}_{\{r\}}=\widetilde{z}_r
\quad \mbox{ for }r\in\set{n},\qquad
\widetilde{p}_{\set{n}}=\prod_{r=1}^n\widetilde{z}_r.
\end{gather*}
For a set $A$, let $\bbone_{A}(x)=1$, when $x\in A$, and 
$\bbone_{A}(x)=0$ otherwise. 
We always set $0^0=1$, $\frac{1}{0}=\infty$, and $1^\infty=1$.
As usual, empty sums, resp.\ empty products, are defined to be zero, 
resp.\ one. 

%%%%%%%%%%%%%%%%%%%%%%%%%%%%%%%%%%%%%%%%%%%%%%%%%%%%%%%%%%%%%%%%%%%%%%
\section{Some identities for permanents}
%%%%%%%%%%%%%%%%%%%%%%%%%%%%%%%%%%%%%%%%%%%%%%%%%%%%%%%%%%%%%%%%%%%%%%
Our first main result is Theorem \ref{t6296} below, 
the proof of which requires the following lemma.
For $r,s\in\set{N}$, let $\tau_{r,s}\in\set{N}_{\neq}^{N}$ be the 
transposition, which interchanges $r$ with $s$, i.e.\ 
\begin{align*}
\tau_{r,s}(\ell)=
\left\{\begin{array}{ll}
\ell&\mbox{ for }\ell\in\set{N}\setminus\{r,s\},\\
s&\mbox{ for }\ell=r,\\
r&\mbox{ for }\ell=s.
\end{array}\right.
\end{align*}
%%%%%%%%%%%%%%%%%%%%%%%%%%%%%%%%%%%%%%%%%%%%%%%%%%%%%%%%%%%%%%%%%%%%%%
\begin{lemma}\label{l437865}
Let $r\in\set{N}$, $T_r=(T_{r,1},T_{r,2}):\,(\set{N}_{\neq}^N)^2
\longrightarrow (\set{N}_{\neq}^N)^2$, 
$T_r(j,k)=(T_{r,1}(j,k),T_{r,2}(j,k))$
with $T_{r,1}(j,k)=j\circ\tau_{r,(j^{-1}\circ k)(r)}$,
$T_{r,2}(j,k)=k\circ\tau_{r,(k^{-1}\circ j)(r)}$ for 
$(j,k)\in(\set{N}_{\neq}^N)^2$. 
Here $\circ$ means composition of functions. 
Then $T_r\circ T_r$ is the identity map on $(\set{N}_{\neq}^N)^2$. 
In particular, $T_r$ is bijective and we have 
$T_{r,1}(j,k)(r)=k_r$, $T_{r,2}(j,k)(r)=j_r$. 
\end{lemma}
%%%%%%%%%%%%%%%%%%%%%%%%%%%%%%%%%%%%%%%%%%%%%%%%%%%%%%%%%%%%%%%%%%%%%%
\Proof 
For $(j,k)\in(\set{N}_{\neq}^N)^2$, set 
$\widetilde{j}=T_{r,1}(j,k)$ and $\widetilde{k}=T_{r,2}(j,k)$. Then
\begin{align*}
(\widetilde{j}^{-1}\circ \widetilde{k})(r)
&=((j\circ\tau_{r,(j^{-1}\circ k)(r)})^{-1}\circ 
k\circ\tau_{r,(k^{-1}\circ j)(r)})(r)\\
&=(\tau_{r,(j^{-1}\circ k)(r)}^{-1}\circ 
j^{-1}\circ k)((k^{-1}\circ j)(r))
=\tau_{r,(j^{-1}\circ k)(r)}(r)
=(j^{-1}\circ k)(r)
\end{align*}
and therefore 
\begin{align*}
T_{r,1}(T_r(j,k))
&=\widetilde{j}\circ\tau_{r,(\widetilde{j}^{-1}
\circ \widetilde{k})(r)}
=(j\circ\tau_{r,(j^{-1}\circ k)(r)})\circ
\tau_{r,(j^{-1}\circ k)(r)}=j.
\end{align*}
Similarly 
$(\widetilde{k}^{-1}\circ \widetilde{j})(r)=(k^{-1}\circ j)(r)$
and $T_{r,2}(T_r(j,k))=k$. \hfill\qed
%%%%%%%%%%%%%%%%%%%%%%%%%%%%%%%%%%%%%%%%%%%%%%%%%%%%%%%%%%%%%%%%%%%%%%
\begin{remark}
\begin{enumerate}

\item Another way of describing $T_r$ is the following: 
For $j,k\in\set{N}_{\neq}^N$, we obtain $T_{r,1}(j,k)$ and 
$T_{r,2}(j,k)$, if in both tuples $j$ and $k$, 
we replace $j_r$ with $k_r$. More precisely, 
if $r\in\set{N}$, $j,k\in\set{N}_{\neq}^N$, 
$a,b\in\set{N}$, $u,v\in\set{N}$ with $j_r=u$, $j_a=v$, $k_r=v$, 
$k_b=u$, then
\begin{align*}
T_{r,1}(j,k)(s) 
=\left\{\begin{array}{ll}
j_s & \mbox{ for }s\in\set{N}\setminus\{r,a\},\\
v & \mbox{ for } s=r,\\
u & \mbox{ for } s=a,
\end{array}\right.
\qquad 
T_{r,2}(j,k)(s) 
=\left\{\begin{array}{ll}
k_s & \mbox{ for }s\in\set{N}\setminus\{r,b\},\\
u & \mbox{ for } s=r,\\
v & \mbox{ for } s=b.
\end{array}\right.
\end{align*}
For example, if $N=3$, $r=2$, $j=(j_1,j_2,j_3)=(2,1,3)$ and
$k=(k_1,k_2,k_3)=(3,2,1)$, then $T_{r,1}(j,k)=(1,2,3)$ and
$T_{r,2}(j,k)=(3,1,2)$. 

\item According to Lemma \ref{l437865}, we have
$(\set{N}_{\neq}^N)^2=\{T_r(j,k)\,|\,(j,k)\in(\set{N}_{\neq}^N)^2\}$
for all $r\in\set{N}$. Hence, for an arbitrary function  
$f:\,(\set{N}_{\neq}^N)^2\longrightarrow\calz$,  
\begin{align*}
\sum_{j\in\set{N}_{\neq}^N}\sum_{k\in\set{N}_{\neq}^N}f(j,k)
=\sum_{j\in\set{N}_{\neq}^N}\sum_{k\in\set{N}_{\neq}^N}f(T_r(j,k)),
\end{align*}
which is the main idea in the proof of the next theorem.
\end{enumerate}
\end{remark}
%%%%%%%%%%%%%%%%%%%%%%%%%%%%%%%%%%%%%%%%%%%%%%%%%%%%%%%%%%%%%%%%%%%%%%
\begin{theorem}\label{t6296}
Let $R,S\subseteq\set{n}$ and
$r\in \set{n}\setminus(R\cup S)\neq\emptyset$. Then
\begin{align}
\lefteqn{\overline{p}_{R\cup\{r\}}\overline{p}_S
-\overline{p}_{R}\overline{p}_{S\cup\{r\}}}\nonumber\\
&=\frac{1}{2}\sum_{s\in\set{n}\setminus\{r\}}
\sum_{(u,v)\in\set{N}_{\neq}^2}y_{u,v,r}y_{u,v,s}
\sum_{\newatop{j\in\set{N}_{\neq}^n}{j_r=u,j_s=v}}
\sum_{\newatop{k\in\set{N}_{\neq}^n}{k_r=u}}
\bigl(\bbone_{S}(s)p_{j,S\setminus\{s\}}p_{k,R}
-\bbone_{R}(s)p_{j,R\setminus\{s\}}p_{k,S}\bigr).\label{e74675}
\end{align}
If $n=1$, then the right-hand side of the equality in 
\eqref{e74675} is defined to be zero. 
\end{theorem}
%%%%%%%%%%%%%%%%%%%%%%%%%%%%%%%%%%%%%%%%%%%%%%%%%%%%%%%%%%%%%%%%%%%%%%
\Proof We have
\begin{align}
\lefteqn{((N-n)!)^2(\overline{p}_{R\cup\{r\}}\overline{p}_S
-\overline{p}_{R}\overline{p}_{S\cup\{r\}})}\nonumber\\ 
&=((N-n)!)^2\sum_{j\in\set{N}_{\neq}^n}\sum_{k\in\set{N}_{\neq}^n}
(p_{j,R\cup\{r\}}p_{k,S}-p_{j,R}p_{k,S\cup\{r\}})\nonumber\\
&=\sum_{j\in\set{N}_{\neq}^N}\sum_{k\in\set{N}_{\neq}^N}
y_{j_r,k_r,r}p_{j,R}p_{k,S}\label{e621950}\\
&=\sum_{a\in\set{N}\setminus\{r\}}\sum_{b\in\set{N}\setminus\{r\}}
\sum_{(u,v)\in\set{N}_{\neq}^2}
\sum_{\newatop{j\in\set{N}_{\neq}^N}{j_r=u,j_a=v}}
\sum_{\newatop{k\in\set{N}_{\neq}^N}{k_r=v,k_b=u}}
y_{u,v,r}p_{j,R}p_{k,S}.\nonumber
\end{align}
Now we use the decompositions 
$\set{N}\setminus\{r\}=(\set{N}\setminus(R\cup\{r\}))\cup R
=(\set{N}\setminus(S\cup\{r\}))\cup S$ 
and obtain 
\begin{align}\label{e46496}
((N-n)!)^2(\overline{p}_{R\cup\{r\}}\overline{p}_S
-\overline{p}_{R}\overline{p}_{S\cup\{r\}})
=A_1+A_2+A_3+A_4,
\end{align}
where
\begin{align*}
A_1&=\sum_{a\in R}\sum_{b\in\set{N}\setminus (S\cup\{r\})}
\sum_{(u,v)\in\set{N}_{\neq}^2}
\sum_{\newatop{j\in\set{N}_{\neq}^N}{j_r=u,j_a=v}}
\sum_{\newatop{k\in\set{N}_{\neq}^N}{k_r=v,k_b=u}}
y_{u,v,r}z_{v,a}p_{j,R\setminus\{a\}}p_{k,S},\\
A_2&=\sum_{a\in \set{N}\setminus (R\cup\{r\})}\sum_{b\in S}
\sum_{(u,v)\in\set{N}_{\neq}^2}
\sum_{\newatop{j\in\set{N}_{\neq}^N}{j_r=u,j_a=v}}
\sum_{\newatop{k\in\set{N}_{\neq}^N}{k_r=v,k_b=u}}
y_{u,v,r}z_{u,b}p_{j,R}p_{k,S\setminus\{b\}},\\
A_3&=\sum_{a\in R}\sum_{b\in S}
\sum_{(u,v)\in\set{N}_{\neq}^2}
\sum_{\newatop{j\in\set{N}_{\neq}^N}{j_r=u,j_a=v}}
\sum_{\newatop{k\in\set{N}_{\neq}^N}{k_r=v,k_b=u}}
y_{u,v,r}z_{v,a}z_{u,b}p_{j,R\setminus\{a\}}
p_{k,S\setminus\{b\}},\\
A_4&=\sum_{a\in \set{N}\setminus (R\cup\{r\})}
\sum_{b\in\set{N}\setminus (S\cup\{r\})}
\sum_{(u,v)\in\set{N}_{\neq}^2}
\sum_{\newatop{j\in\set{N}_{\neq}^N}{j_r=u,j_a=v}}
\sum_{\newatop{k\in\set{N}_{\neq}^N}{k_r=v,k_b=u}}
y_{u,v,r}p_{j,R}p_{k,S}.
\end{align*}
A representation similar to \eqref{e46496} can be shown by using
\eqref{e621950}, Lemma \ref{l437865} and the fact that 
$y_{v,u,r}=-y_{u,v,r}$ for $u,v\in\set{N}$. Indeed, if
$T_r=(T_{r,1},T_{r,2})$ is defined as in that lemma, then
\begin{align}
\lefteqn{((N-n)!)^2(\overline{p}_{R\cup\{r\}}\overline{p}_S
-\overline{p}_{R}\overline{p}_{S\cup\{r\}})
=\sum_{j\in\set{N}_{\neq}^N}\sum_{k\in\set{N}_{\neq}^N}
y_{k_r,j_r,r}p_{T_{r,1}(j,k),R}p_{T_{r,2}(j,k),S}
}\nonumber\\ 
&\hspace{3cm}
=-\sum_{a\in\set{N}\setminus\{r\}}\sum_{b\in\set{N}\setminus\{r\}}
\sum_{(u,v)\in\set{N}_{\neq}^2}
\sum_{\newatop{j\in\set{N}_{\neq}^N}{j_r=u,j_a=v}}
\sum_{\newatop{k\in\set{N}_{\neq}^N}{k_r=v,k_b=u}}
y_{u,v,r} p_{T_{r,1}(j,k),R}p_{T_{r,2}(j,k),S}\nonumber\\
&\hspace{3cm}
=-(A_1'+A_2'+A_3'+A_4'),\label{e729587}
\end{align}
where 
\begin{align*}
A_1'&=\sum_{a\in R}\sum_{b\in\set{N}\setminus (S\cup\{r\})}
\sum_{(u,v)\in\set{N}_{\neq}^2}
\sum_{\newatop{j\in\set{N}_{\neq}^N}{j_r=u,j_a=v}}
\sum_{\newatop{k\in\set{N}_{\neq}^N}{k_r=v,k_b=u}}
y_{u,v,r}z_{u,a}p_{j,R\setminus\{a\}}p_{k,S},\\
A_2'&=\sum_{a\in \set{N}\setminus (R\cup\{r\})}\sum_{b\in S}
\sum_{(u,v)\in\set{N}_{\neq}^2}
\sum_{\newatop{j\in\set{N}_{\neq}^N}{j_r=u,j_a=v}}
\sum_{\newatop{k\in\set{N}_{\neq}^N}{k_r=v,k_b=u}}
y_{u,v,r}z_{v,b}p_{j,R}p_{k,S\setminus\{b\}},\\
A_3'&=\sum_{a\in R}\sum_{b\in S}\sum_{(u,v)\in\set{N}_{\neq}^2}
\sum_{\newatop{j\in\set{N}_{\neq}^N}{j_r=u,j_a=v}}
\sum_{\newatop{k\in\set{N}_{\neq}^N}{k_r=v,k_b=u}}
y_{u,v,r}z_{u,a}z_{v,b}p_{j,R\setminus\{a\}}p_{k,S\setminus\{b\}}
\end{align*}
and $A_4'=A_4$. 
Adding the right-hand sides of \eqref{e46496} and 
\eqref{e729587} and dividing by two, we get the identity
\begin{align}\label{e6438697}
((N-n)!)^2(\overline{p}_{R\cup\{r\}}\overline{p}_S
-\overline{p}_{R}\overline{p}_{S\cup\{r\}})
&=\frac{1}{2}(-B_1+B_2+B_3),
\end{align}
where
\begin{align*}
B_1&=-A_1+A_1'= 
\sum_{a\in R}\sum_{b\in\set{N}\setminus (S\cup\{r\})}
\sum_{(u,v)\in\set{N}_{\neq}^2}
\sum_{\newatop{j\in\set{N}_{\neq}^N}{j_r=u,j_a=v}}
\sum_{\newatop{k\in\set{N}_{\neq}^N}{k_r=v,k_b=u}}
y_{u,v,r}y_{u,v,a}p_{j,R\setminus\{a\}}p_{k,S},\\
B_2&=A_2-A_2'= 
\sum_{a\in \set{N}\setminus (R\cup\{r\})}\sum_{b\in S}
\sum_{(u,v)\in\set{N}_{\neq}^2}
\sum_{\newatop{j\in\set{N}_{\neq}^N}{j_r=u,j_a=v}}
\sum_{\newatop{k\in\set{N}_{\neq}^N}{k_r=v,k_b=u}}
y_{u,v,r}y_{u,v,b}p_{j,R}p_{k,S\setminus\{b\}},\\
B_3&=A_3-A_3'=\sum_{a\in R}\sum_{b\in S}
\sum_{(u,v)\in\set{N}_{\neq}^2}
\sum_{\newatop{j\in\set{N}_{\neq}^N}{j_r=u,j_a=v}}
\sum_{\newatop{k\in\set{N}_{\neq}^N}{k_r=v,k_b=u}}
y_{u,v,r}(z_{v,a}z_{u,b}-z_{u,a}z_{v,b})p_{j,R\setminus\{a\}}
p_{k,S\setminus\{b\}}. 
\end{align*}
Now we write $B_3=B_3'-B_3''$, where
\begin{align}
B_3'&=\sum_{a\in R}\sum_{b\in S}
\sum_{(u,v)\in\set{N}_{\neq}^2}
\sum_{\newatop{j\in\set{N}_{\neq}^N}{j_r=u,j_a=v}}
\sum_{\newatop{k\in\set{N}_{\neq}^N}{k_r=v,k_b=u}}
y_{u,v,r}z_{v,a}y_{u,v,b}p_{j,R\setminus\{a\}}p_{k,S\setminus\{b\}}
\nonumber\\
&=\sum_{a\in R}\sum_{b\in S}
\sum_{(u,v)\in\set{N}_{\neq}^2}
\sum_{\newatop{j\in\set{N}_{\neq}^N}{j_r=u,j_a=v}}
\sum_{\newatop{k\in\set{N}_{\neq}^N}{k_r=v,k_b=u}}
y_{u,v,r}y_{u,v,b}p_{j,R}p_{k,S\setminus\{b\}}\label{e73289650}
\end{align} 
and
\begin{align}
B_3''
&=\sum_{a\in R}\sum_{b\in S}
\sum_{(u,v)\in\set{N}_{\neq}^2}
\sum_{\newatop{j\in\set{N}_{\neq}^N}{j_r=u,j_a=v}}
\sum_{\newatop{k\in\set{N}_{\neq}^N}{k_r=v,k_b=u}}
y_{u,v,r}y_{u,v,a}z_{v,b}p_{j,R\setminus\{a\}}p_{k,S\setminus\{b\}}
\label{e53760}\\
&=\sum_{a\in R}\sum_{b\in S}
\sum_{(u,v)\in\set{N}_{\neq}^2}
\sum_{\newatop{j\in\set{N}_{\neq}^N}{j_r=u,j_a=v}}
\sum_{\newatop{k\in\set{N}_{\neq}^N}{k_r=v,k_b=u}}
y_{u,v,r}y_{u,v,a}p_{j,R\setminus\{a\}}p_{k,S}.\label{e756860}
\end{align}
Indeed, \eqref{e756860} can be derived from \eqref{e53760}
by interchanging $u$ with $v$ for 
$(a,b)\in R\times S$ being fixed. We note that here 
$y_{u,v,r}y_{u,v,a}=y_{v,u,r}y_{v,u,a}$ and 
\begin{align*}
\sum_{\newatop{j\in\set{N}_{\neq}^N}{j_r=v,j_a=u}}
p_{j,R\setminus\{a\}}
&=\sum_{\newatop{j\in\set{N}_{\neq}^N}{j_r=u,j_a=v}}
p_{j,R\setminus\{a\}},\qquad
\sum_{\newatop{k\in\set{N}_{\neq}^N}{k_r=u,k_b=v}}z_{u,b}
p_{k,S\setminus\{b\}}
=\sum_{\newatop{k\in\set{N}_{\neq}^N}{k_r=v,k_b=u}}p_{k,S}
\end{align*}
for $(u,v)\in\set{N}_{\neq}^2$.
Combining \eqref{e6438697}, \eqref{e73289650} and \eqref{e756860},
we get 
\begin{align}\label{e8166}
((N-n)!)^2(\overline{p}_{R\cup\{r\}}\overline{p}_S
-\overline{p}_{R}\overline{p}_{S\cup\{r\}})=\frac{1}{2}(C_1-C_2),
\end{align}
where $C_1=B_2+B_3'$, $C_2=B_1+B_3''$. We have 
\begin{align}
C_1&=\sum_{a\in \set{N}\setminus\{r\}}\sum_{b\in S}
\sum_{(u,v)\in\set{N}_{\neq}^2}
\sum_{\newatop{j\in\set{N}_{\neq}^N}{j_r=u,j_a=v}}
\sum_{\newatop{k\in\set{N}_{\neq}^N}{k_r=v,k_b=u}}
y_{u,v,r}y_{u,v,b}p_{j,R}p_{k,S\setminus\{b\}}\label{e626559}\\
&=\sum_{s\in\set{n}\setminus \{r\}}
\sum_{(u,v)\in\set{N}_{\neq}^2}y_{u,v,r}y_{u,v,s}
\sum_{\newatop{j\in\set{N}_{\neq}^N}{j_r=u}}
\sum_{\newatop{k\in\set{N}_{\neq}^N}{k_r=v,k_s=u}}
\bbone_{S}(s)p_{j,R}p_{k,S\setminus\{s\}}\label{e214765}\\
&=\sum_{s\in\set{n}\setminus \{r\}}
\sum_{(u,v)\in\set{N}_{\neq}^2}y_{u,v,r}y_{u,v,s}
\sum_{\newatop{j\in\set{N}_{\neq}^N}{j_r=v,j_s=u}}
\sum_{\newatop{k\in\set{N}_{\neq}^N}{k_r=u}}
\bbone_{S}(s)p_{j,S\setminus\{s\}}p_{k,R}\label{e936659}\\
&=\sum_{s\in\set{n}\setminus \{r\}}
\sum_{(u,v)\in\set{N}_{\neq}^2}y_{u,v,r}y_{u,v,s}
\sum_{\newatop{j\in\set{N}_{\neq}^N}{j_r=u,j_s=v}}
\sum_{\newatop{k\in\set{N}_{\neq}^N}{k_r=u}}
\bbone_{S}(s)p_{j,S\setminus\{s\}}p_{k,R}.\label{e4321353}
\end{align}
Here, \eqref{e626559} follows from the definitions of $B_2$ 
and $B_3'$. To get \eqref{e214765}, we replace $b$ by $s$
and note that 
$\sum_{a\in \set{N}\setminus\{r\}} 
\sum_{{j\in\set{N}_{\neq}^N}:\,{j_r=u,j_a=v}}
=\sum_{{j\in\set{N}_{\neq}^N}:\,{j_r=u}}$ 
for fixed $(u,v)\in\set{N}_{\neq}^2$. 
For \eqref{e936659}, we interchanged $j$ with~$k$.
Finally, \eqref{e4321353} follows by interchanging $j_r$ with $j_s$
and noting that 
$\sum_{{j\in\set{N}_{\neq}^N}:\,{j_r=v,j_s=u}}
p_{j,S\setminus\{s\}}
=\sum_{{j\in\set{N}_{\neq}^N}:\,{j_r=u,j_s=v}}
p_{j,S\setminus\{s\}}$ for fixed $(u,v)\in\set{N}_{\neq}^2$.
Similarly, 
\begin{align}
C_2&=\sum_{a\in R}\sum_{b\in\set{N}\setminus\{r\}}
\sum_{(u,v)\in\set{N}_{\neq}^2}
\sum_{\newatop{j\in\set{N}_{\neq}^N}{j_r=u,j_a=v}}
\sum_{\newatop{k\in\set{N}_{\neq}^N}{k_r=v,k_b=u}}
y_{u,v,r}y_{u,v,a}p_{j,R\setminus\{a\}}p_{k,S}\nonumber\\
&=\sum_{s\in \set{n}\setminus\{r\}}
\sum_{(u,v)\in\set{N}_{\neq}^2}y_{u,v,r}y_{u,v,s}
\sum_{\newatop{j\in\set{N}_{\neq}^N}{j_r=u,j_s=v}}
\sum_{\newatop{k\in\set{N}_{\neq}^N}{k_r=v}}
\bbone_{R}(s)p_{j,R\setminus\{s\}}p_{k,S}\nonumber\\
&=\sum_{s\in \set{n}\setminus\{r\}}
\sum_{(u,v)\in\set{N}_{\neq}^2}y_{u,v,r}y_{u,v,s}
\sum_{\newatop{j\in\set{N}_{\neq}^N}{j_r=v,j_s=u}}
\sum_{\newatop{k\in\set{N}_{\neq}^N}{k_r=u}}
\bbone_{R}(s)p_{j,R\setminus\{s\}}p_{k,S}\label{e3658}\\
&=\sum_{s\in \set{n}\setminus\{r\}}
\sum_{(u,v)\in\set{N}_{\neq}^2}y_{u,v,r}y_{u,v,s}
\sum_{\newatop{j\in\set{N}_{\neq}^N}{j_r=u,j_s=v}}
\sum_{\newatop{k\in\set{N}_{\neq}^N}{k_r=u}}
\bbone_{R}(s)p_{j,R\setminus\{s\}}p_{k,S}.\label{e385025}
\end{align}
In particular, \eqref{e3658} follows by interchanging $u$ with $v$
for fixed $s\in\set{n}\setminus\{r\}$.
Combining \eqref{e8166}, \eqref{e4321353}, and \eqref{e385025} the 
assertion is shown.\hfill\qed  \medskip
%%%%%%%%%%%%%%%%%%%%%%%%%%%%%%%%%%%%%%%%%%%%%%%%%%%%%%%%%%%%%%%%%%%%%%

The next result on elementary symmetric polynomials
is due to Dougall \citet[formula~(3) on page~65]{Dougall1905}. 
We now show that it is a consequence of Theorem \ref{t6296}. 
%%%%%%%%%%%%%%%%%%%%%%%%%%%%%%%%%%%%%%%%%%%%%%%%%%%%%%%%%%%%%%%%%%%%%%

\begin{corollary}\label{c76276597}
Let $N\in\NN$, $z_j\in\calz$ for $j\in\set{N}$, 
$E_{A,k}=\sum_{\newatop{J\subseteq A}{\card{J}=k}}\prod_{j\in J}z_j$
for $A\subseteq\set{N}$ and $k\in\ZZ$, 
$a,b\in\setn{N}=\{0,1,\dots,N\}$. In particular, $E_{A,0}=1$ and
$E_{A,k}=0$ if $k<0$ or $k>\card{A}$. Then 
\begin{align*}
\lefteqn{
(a+1)(N-b)E_{\set{N},a+1}E_{\set{N},b}
-(b+1)(N-a)E_{\set{N},a}E_{\set{N},b+1}}\\
&\hspace{3cm}=\frac{1}{2}
\sum_{(u,v)\in\set{N}_{\neq}^2}(z_u-z_v)^2
(E_{\set{N}\setminus\{u,v\},b-1}E_{\set{N}\setminus\{u,v\},a}
-E_{\set{N}\setminus\{u,v\},a-1}E_{\set{N}\setminus\{u,v\},b}).
\end{align*}
\end{corollary}
%%%%%%%%%%%%%%%%%%%%%%%%%%%%%%%%%%%%%%%%%%%%%%%%%%%%%%%%%%%%%%%%%%%%%
\Proof
For $a=N$ or $b=N$, the assertion is trivial. Let now
$a<N$ and $b<N$ and consider the assumptions of Theorem \ref{t6296}, 
where $n=N$ and
$Z$ has identical columns, i.e.\ $z_{j,1}=\dots=z_{j,n}=z_j$
for all $j\in\set{N}$. Further, let $\card{R}=a$ and $\card{S}=b$. 
Then the assertion follows from \eqref{e74675} and 
\begin{gather*}
\overline{p}_{R}
=a!(N-a)!E_{\set{N},a},\quad
\overline{p}_{S}
=b!(N-b)!E_{\set{N},b},\\
\overline{p}_{R\cup\{r\}}
=(a+1)!(N-a-1)!E_{\set{N},a+1},\quad
\overline{p}_{S\cup\{r\}}
=(b+1)!(N-b-1)!E_{\set{N},b+1},\\
\sum_{\ell\in S}\sum_{\newatop{j\in\set{N}_{\neq}^n}{j_r=u,j_\ell=v}}
p_{j,S\setminus\{\ell\}}
=b!(N-b-1)!E_{\set{N}\setminus\{u,v\},b-1},\quad
\sum_{\newatop{k\in\set{N}_{\neq}^n}{k_r=u}}p_{k,R}
=a!(N-a-1)!E_{\set{N}\setminus\{u\},a},\\
\sum_{\ell\in R}\sum_{\newatop{j\in\set{N}_{\neq}^n}{j_r=u,j_\ell=v}}
p_{j,R\setminus\{\ell\}}
=a!(N-a-1)!E_{\set{N}\setminus\{u,v\},a-1},\quad
\sum_{\newatop{k\in\set{N}_{\neq}^n}{k_r=u}}p_{k,S}
=b!(N-b-1)!E_{\set{N}\setminus\{u\},b}
\end{gather*}
and $E_{\set{N}\setminus\{u\},a}
=E_{\set{N}\setminus\{u,v\},a}+z_vE_{\set{N}\setminus\{u,v\},a-1}$,
$E_{\set{N}\setminus\{u\},b}
=E_{\set{N}\setminus\{u,v\},b}+z_vE_{\set{N}\setminus\{u,v\},b-1}$,
where $(u,v)\in\set{N}_{\neq}^2$.\hfill\qed

%%%%%%%%%%%%%%%%%%%%%%%%%%%%%%%%%%%%%%%%%%%%%%%%%%%%%%%%%%%%%%%%%%%%%%
\begin{corollary}\label{c42875}
Let $R\subseteq\set{n}$ and $r\in \set{n}\setminus R\neq \emptyset$. 
Then  
\begin{align}\label{e52674}
\overline{p}_{R\cup\{r\}}-\widetilde{z}_r\overline{p}_{R}
=-\frac{1}{2N}\sum_{s\in R}\sum_{j\in\set{N}_{\neq}^{n}}
y_{j_r,j_s,r}y_{j_r,j_s,s}p_{j,R\setminus\{s\}}.
\end{align}
\end{corollary}
%%%%%%%%%%%%%%%%%%%%%%%%%%%%%%%%%%%%%%%%%%%%%%%%%%%%%%%%%%%%%%%%%%%%%%
\Proof In Theorem \ref{t6296}, set $S=\emptyset$. \hfill\qed
\medskip

%%%%%%%%%%%%%%%%%%%%%%%%%%%%%%%%%%%%%%%%%%%%%%%%%%%%%%%%%%%%%%%%%%%%%%
The next result follows from Corollary \ref{c42875}
and is the main argument in the proof of our inequalities 
in Section \ref{s4518u5}.

%%%%%%%%%%%%%%%%%%%%%%%%%%%%%%%%%%%%%%%%%%%%%%%%%%%%%%%%%%%%%%%%%%%%%%
\begin{theorem}\label{t857449}
If $r\in\set{n}_{\neq}^n$, 
$R_{k}:=R_k(r):=\{r_\ell\,|\,\ell\in\set{k}\}$
for $k\in\setn{n}=\{0,\dots,n\}$, then $\emptyset=R_0\subsetneq R_1
\subsetneq\dots\subsetneq R_n=\set{n}$ 
is a maximal chain of subsets of $\set{n}$ and
\begin{align}
\overline{p}_{\set{n}}
-\frac{N!\,\widetilde{p}_{\set{n}}}{(N-n)!}
&=-\frac{1}{2N}\sum_{k=2}^n
\sum_{s\in R_{k-1}}\sum_{j\in\set{N}_{\neq}^{n}}
y_{j_{r_k},j_s,r_k}y_{j_{r_k},j_s,s}p_{j,R_{k-1}\setminus\{s\}}
\widetilde{p}_{\set{n}\setminus R_k}. \label{e783564}
\end{align}
On the other hand, 
\begin{align}
\lefteqn{\overline{p}_{\set{n}}
-\frac{N!\,\widetilde{p}_{\set{n}}}{(N-n)!}
=-\sum_{k=2}^n\frac{1}{2Nk\binomial{n}{k}}
\sum_{\newatop{R\subseteq\set{n}}{\card{R}=k}}
\sum_{(r,s)\in R_{\neq}^2}\sum_{j\in\set{N}_{\neq}^{n}}
y_{j_r,j_{s},r}y_{j_r,j_{s},s}p_{j,R\setminus\{r,s\}}
\widetilde{p}_{\set{n}\setminus R}}\label{e8366548}\\
&=-\sum_{(u,v)\in\set{N}_{\neq}^2}\sum_{(r,s)\in\set{n}_{\neq}^2}
y_{u,v,r}y_{u,v,s}
\sum_{k=2}^n\frac{1}{2Nk\binomial{n}{k}}
\sum_{\newatop{R\subseteq\set{n}\setminus\{r,s\}}{\card{R}=k-2}}
\widetilde{p}_{\set{n}\setminus (R\cup\{r,s\})}
\sum_{j\in(\set{N}\setminus\{u,v\})_{\neq}^{
\set{n}\setminus\{r,s\}}}p_{j,R}.\label{e5703}
\end{align}
If $n=1$, the right-hand sides of \eqref{e783564},
\eqref{e8366548} and \eqref{e5703} are 
defined to be zero. 
\end{theorem}
%%%%%%%%%%%%%%%%%%%%%%%%%%%%%%%%%%%%%%%%%%%%%%%%%%%%%%%%%%%%%%%%%%%%%%
\Proof 
In view of Corollary \ref{c42875} and the identities 
$\overline{p}_{\emptyset}=\frac{N!}{(N-n)!}$ and
$\widetilde{p}_{\emptyset}=1$, we see that
\begin{align*}
\overline{p}_{\set{n}}-
\frac{N!\,\widetilde{p}_{\set{n}}}{(N-n)!}
&=\sum_{k=1}^n
(\widetilde{p}_{\set{n}\setminus R_k}\overline{p}_{R_k}
-\widetilde{p}_{\set{n}\setminus R_{k-1}}\overline{p}_{R_{k-1}})\\
&=\sum_{k=1}^n(\overline{p}_{R_{k}}
-\widetilde{z}_{r_k}\overline{p}_{R_{k-1}})
\widetilde{p}_{\set{n}\setminus R_k}\\
&=-\frac{1}{2N}\sum_{k=2}^n
\sum_{s\in R_{k-1}}\sum_{j\in\set{N}_{\neq}^{n}}
y_{j_{r_k},j_s,r_k}y_{j_{r_k},j_s,s}p_{j,R_{k-1}\setminus\{s\}}
\widetilde{p}_{\set{n}\setminus R_k},
\end{align*}
giving \eqref{e783564}. Hence,
\begin{align*}
\overline{p}_{\set{n}}-\frac{N!\,\widetilde{p}_{\set{n}}}{(N-n)!}
&=-\frac{1}{2N}\sum_{k=2}^n\frac{1}{n!}
\sum_{r\in\set{n}_{\neq}^n}
\sum_{s\in R_{k}(r)\setminus\{r_k\}}\sum_{j\in\set{N}_{\neq}^{n}}
y_{j_{r_k},j_s,r_k}y_{j_{r_k},j_s,s}p_{j,R_k(r)\setminus\{s,r_k\}}
\widetilde{p}_{\set{n}\setminus R_k(r)}\\
&=-\sum_{k=2}^n\frac{1}{2Nk\binomial{n}{k}}
\sum_{\newatop{R\subseteq\set{n}}{\card{R}=k}}\sum_{t\in R}
\sum_{s\in R\setminus\{t\}}\sum_{j\in\set{N}_{\neq}^{n}}
y_{j_t,j_s,t}y_{j_t,j_s,s}p_{j,R\setminus\{s,t\}}
\widetilde{p}_{\set{n}\setminus R},
\end{align*}
since for $k\in\set{n}\setminus\{1\}$,
$R\subseteq n$ with $\card{R}=k$ and $t\in R$, the number of 
$r\in\set{n}_{\neq}^n$ with $R_k(r)=R$ and $r_k=t$
is equal to $(k-1)!(n-k)!$. This shows \eqref{e8366548}. 
Furthermore, \eqref{e5703} is clear. \hfill\qed
%%%%%%%%%%%%%%%%%%%%%%%%%%%%%%%%%%%%%%%%%%%%%%%%%%%%%%%%%%%%%%%%%%%%%%

\begin{remark} 
\begin{enumerate}

\item If $Z$ has identical rows, i.e.\ $z_{1,r}=\dots=z_{N,r}$
for all $r\in\set{n}$, then 
$\overline{p}_{R}=\frac{N!\widetilde{p}_{R}}{(N-n)!}$
for all $R\subseteq\set{n}$ and therefore both sides in each 
identity \eqref{e74675} and \eqref{e52674}--\eqref{e5703} give zero. 
For identities in the case, when $Z$ has identical columns, 
see Corollary \ref{c526745} below.

\item The identities of Theorem \ref{t857449} can be rewritten 
as expansions for the permanent $\overline{p}_{\set{n}}$. 
Further such formulas can be found in the literature, 
e.g.\ see Minc \citet[Chapter 7]{MR504978}. For instance, 
Ryser \citet[Theorem 4.1, page 26]{MR0150048} proved that 
\begin{align*}
\overline{p}_{\set{n}}=\sum_{k=1}^n(-1)^{n-k}
\binomial{N-k}{n-k}
\sum_{\newatop{J\subseteq \set{N}}{\card{J}=k}}
\prod_{r=1}^n\Big(\sum_{j\in J}z_{j,r}\Big).
\end{align*}
In the case $n=N$, this implies that
$\overline{p}_{\set{n}}-n^n\widetilde{p}_{\set{n}}
=\sum_{J\subseteq \set{n}:\,1\leq\card{J}<n}(-1)^{n-\card{J}}
\prod_{r=1}^n\sum_{j\in J}z_{j,r}$, which however 
is not comparable with the identities of Theorem \ref{t857449}
under the present assumption.
We note that a second order expansion for $\overline{p}_{\set{n}}$
can be found in Theorem \ref{t384696} below.

\item 
Let us assume that $Z=(z_{j,r})\in[0,\infty)^{N\times n}$
has decreasing columns, i.e.\ 
$z_{j,r}\geq z_{j+1,r}$ for all $j\in\set{n-1}$ and $r\in\set{n}$.
Then $y_{j_1,j_2,r}y_{j_1,j_2,s}\geq0$ 
for all $j_1,j_2\in\set{N}$ and $r,s\in\set{n}$. 
Therefore, Corollary \ref{c42875} implies in this case that
$\overline{p}_{R\cup\{r\}}\leq\widetilde{z}_r\overline{p}_{R}$
for $R\subseteq\set{n}$ and $r\in\set{n}\setminus R\neq\emptyset$. 
Further, Theorem \ref{t857449} gives  in this case that
$\overline{p}_{\set{n}}
\leq\frac{N!\,\widetilde{p}_{\set{n}}}{(N-n)!}$.
Both inequalities above also follow from the more general 
Corollary 4.9 in Br\"and\'{e}n et al.\ \citet{MR3051161}, which was 
shown with the help of 
the monotone column permanent theorem. 

\end{enumerate}
\end{remark}
%%%%%%%%%%%%%%%%%%%%%%%%%%%%%%%%%%%%%%%%%%%%%%%%%%%%%%%%%%%%%%%%%%%%%%
 
\begin{corollary} \label{c526745}
Let $N\in\NN$, $n\in\set{N}$, $z_j\in\calz$ for $j\in\set{N}$,
and $\widetilde{z}=\frac{1}{N}\sum_{j=1}^Nz_j$. 
For $A\subseteq\set{N}$ and $k\in\ZZ$, let
$E_{A,k}=\sum_{\newatop{J\subseteq A}{\card{J}=k}}\prod_{j\in J}z_j$.
Then
\begin{align}\label{e141649}
\frac{1}{\binomial{N}{n}}E_{\set{N},n}-\widetilde{z}^n
=-\frac{1}{2N}
\sum_{(u,v)\in\set{N}_{\neq}^2}(z_u-z_v)^2
\sum_{k=2}^n \frac{\widetilde{z}^{n-k}}{k\binomial{N}{k}}
E_{\set{N}\setminus\{u,v\},k-2}.
\end{align}
In particular for $n=N$, we get Dougall's \cite[page 77]{Dougall1905} 
identity 
\begin{align}\label{e51856}
\prod_{j=1}^nz_j-\widetilde{z}^n
&=-\frac{1}{2n}\sum_{(u,v)\in\set{n}_{\neq}^2}(z_u-z_v)^2
\sum_{k=2}^n\frac{\widetilde{z}^{n-k}}{k\binomial{n}{k}}
E_{\set{n}\setminus\{u,v\},k-2}.
\end{align} 
\end{corollary}
%%%%%%%%%%%%%%%%%%%%%%%%%%%%%%%%%%%%%%%%%%%%%%%%%%%%%%%%%%%%%%%%%%%%%%
\Proof 
Identity \eqref{e141649} follows from \eqref{e5703} in the case that
$Z$ has identical columns. Indeed, letting 
$z_{j,1}=\dots=z_{j,n}=z_j$ for all $j\in\set{N}$, then
$\widetilde{z}_1=\dots=\widetilde{z}_n=\widetilde{z}$ and
$\overline{p}_{\set{n}}=n!E_{\set{N},n}$,
$\widetilde{p}_{\set{n}}=\widetilde{z}^n$ and
\begin{align*}
\overline{p}_{\set{n}}
-\frac{N!\,\widetilde{p}_{\set{n}}}{(N-n)!}
&=-\sum_{(u,v)\in\set{N}_{\neq}^2}(z_u-z_v)^2
\sum_{k=2}^n\frac{\widetilde{z}^{n-k}}{2Nk\binomial{n}{k}}
\sum_{(r,s)\in\set{n}_{\neq}^2}
\sum_{\newatop{R\subseteq\set{n}\setminus\{r,s\}}{\card{R}=k-2}}
\sum_{j\in(\set{N}\setminus\{u,v\})_{\neq}^{
\set{n}\setminus\{r,s\}}}\prod_{\ell\in R}z_{j_\ell}\\
&=-\frac{N!}{(N-n)!}\sum_{(u,v)\in\set{N}_{\neq}^2}(z_u-z_v)^2
\sum_{k=2}^n\frac{\widetilde{z}^{n-k}}{2Nk\binomial{N}{k}}
\sum_{\newatop{J\subseteq\set{N}\setminus\{u,v\}}{\card{J}=k-2}}
\prod_{j\in J}z_{j}. 
\end{align*}
Identity \eqref{e51856} follows from \eqref{e141649}, if $n=N$.
\hfill\qed\medskip 
%%%%%%%%%%%%%%%%%%%%%%%%%%%%%%%%%%%%%%%%%%%%%%%%%%%%%%%%%%%%%%%%%%%%%%

\noindent 
We note that the right-hand side of \eqref{e141649}
gives an expansion for the difference 
between the normalized elementary symmetric polynomial
$\frac{1}{\binomial{N}{n}}E_{\set{N},n}$
and $\widetilde{z}^n$.
Further, identities similar to \eqref{e141649} or
\eqref{e51856} have been proved by 
Hurwitz \citet{MR1580244} and Dinghas \citet{MR0023793}.

The next lemma is needed in the proof of our last main result of this 
section. 

%%%%%%%%%%%%%%%%%%%%%%%%%%%%%%%%%%%%%%%%%%%%%%%%%%%%%%%%%%%%%%%%%%%%%%
\begin{lemma}\label{l457638}
If $n\geq 3$, $R\subseteq\set{n}$ with $\card{R}\leq n-3$
and $(r,s,t)\in (\set{n}\setminus R)_{\neq}^3$, then
\begin{align}\label{e6219767}
\sum_{j\in\set{N}_{\neq}^n}y_{j_r,j_s,r}y_{j_r,j_s,s}
(p_{j,R\cup\{t\}}-\widetilde{z}_tp_{j,R})
&=D_1-D_2,
\end{align}
where
\begin{align}
D_1
&:=D_1(r,s,t,R)
:=\frac{2}{N}\sum_{j\in\set{N}_{\neq}^n}
y_{j_r,j_s,r}y_{j_r,j_s,s}y_{j_t,j_r,t}p_{j,R},\label{e186573}\\
D_2
&:=D_2(r,s,t,R)
:=\frac{1}{2N}\sum_{q\in R}\sum_{j\in\set{N}_{\neq}^n}
y_{j_r,j_s,r}y_{j_r,j_s,s}y_{j_t,j_q,t}y_{j_t,j_q,q}
p_{j,R\setminus\{q\}}.\label{e186574}
\end{align}
\end{lemma}
%%%%%%%%%%%%%%%%%%%%%%%%%%%%%%%%%%%%%%%%%%%%%%%%%%%%%%%%%%%%%%%%%%%%%%
\Proof 
Let $D_0:=D_0(r,s,t,R)$ denote the left-hand side of the equation in 
\eqref{e6219767}. For $j\in\set{N}_{\neq}^N$, we have 
$\widetilde{z}_t=\frac{1}{N}\sum_{q=1}^Nz_{j_q,t}$ 
and therefore 
\begin{align*}
(N-n)!D_0 
&=\sum_{j\in\set{N}_{\neq}^N}y_{j_r,j_s,r}y_{j_r,j_s,s}
(z_{j_t,t}-\widetilde{z}_t)p_{j,R}
=\frac{1}{N}\sum_{q\in\set{N}\setminus\{t\}}
\sum_{j\in\set{N}_{\neq}^N}y_{j_r,j_s,r}y_{j_r,j_s,s}
y_{j_t,j_q,t}p_{j,R},
\end{align*}
where 
\begin{align*}
\sum_{q\in\set{N}\setminus\set{n}}
\sum_{j\in\set{N}_{\neq}^N}y_{j_r,j_s,r}y_{j_r,j_s,s}
y_{j_t,j_q,t}p_{j,R}=0,
\end{align*}
which follows by interchanging $j_t$ with $j_q$. Hence
\begin{align*}
D_0
&=\frac{1}{N}\sum_{q\in\set{n}\setminus\{t\}}
\sum_{j\in\set{N}_{\neq}^n}
y_{j_r,j_s,r}y_{j_r,j_s,s}y_{j_t,j_q,t}p_{j,R}
=D_1-D_2+D_3,
\end{align*}
where
\begin{align*}
D_1&=\frac{1}{N}\sum_{q\in\{r,s\}}\sum_{j\in\set{N}_{\neq}^n}
y_{j_r,j_s,r}y_{j_r,j_s,s}y_{j_t,j_q,t}p_{j,R},\\
D_2&=-\frac{1}{N}
\sum_{q\in R}\sum_{j\in\set{N}_{\neq}^n}
y_{j_r,j_s,r}y_{j_r,j_s,s}y_{j_t,j_q,t}p_{j,R},\\
D_3&=\frac{1}{N}
\sum_{q\in \set{n}\setminus (R\cup\{r,s,t\})}
\sum_{j\in\set{N}_{\neq}^n}
y_{j_r,j_s,r}y_{j_r,j_s,s}y_{j_t,j_q,t}p_{j,R}.
\end{align*}
Let us consider the term $D_1$. 
Interchanging $j_s$ with $j_r$ in the summand for $q=s$, we obtain
\begin{align*}
D_1
&=\frac{2}{N}\sum_{j\in\set{N}_{\neq}^n}
y_{j_r,j_s,r}y_{j_r,j_s,s}y_{j_t,j_r,t}p_{j,R}. 
\end{align*}
The term $D_2$ can be treated similarly. By interchanging 
$j_q$ with $j_t$ in the second sum, we derive 
\begin{align}
D_2
&=-\frac{1}{N}
\sum_{q\in R}\sum_{j\in\set{N}_{\neq}^n}
y_{j_r,j_s,r}y_{j_r,j_s,s}y_{j_t,j_q,t}z_{j_q,q}
p_{j,R\setminus\{q\}}\label{e74657}\\
&=\frac{1}{N}
\sum_{q\in R}\sum_{j\in\set{N}_{\neq}^n}
y_{j_r,j_s,r}y_{j_r,j_s,s}y_{j_t,j_q,t}z_{j_t,q}
p_{j,R\setminus\{q\}}.\label{e74658}
\end{align}
Now, adding the right-hand sides of \eqref{e74657}, 
\eqref{e74658} and dividing by two we get
\begin{align*}
D_2
&=\frac{1}{2N}\sum_{q\in R}
\sum_{j\in\set{N}_{\neq}^n}
y_{j_r,j_s,r}y_{j_r,j_s,s}y_{j_t,j_q,t}y_{j_t,j_q,q}
p_{j,R\setminus\{q\}}.
\end{align*}
Finally, we have $D_3=0$, since
\begin{align*}
D_3&=\frac{1}{N}
\sum_{q\in \set{n}\setminus (R\cup\{r,s,t\})}
\sum_{j\in\set{N}_{\neq}^n}
y_{j_r,j_s,r}y_{j_r,j_s,s}y_{j_q,j_t,t}p_{j,R}
=-D_3,
\end{align*}
which follows by interchanging $j_t$ with $j_q$ in the second sum. 
This completes the proof.\hfill\qed\medskip

%%%%%%%%%%%%%%%%%%%%%%%%%%%%%%%%%%%%%%%%%%%%%%%%%%%%%%%%%%%%%%%%%%%%%%
The next result contains a second order expansion for 
$\overline{p}_{\set{n}}$ and is the main argument in the proof of 
Theorem \ref{th487698}. 
%%%%%%%%%%%%%%%%%%%%%%%%%%%%%%%%%%%%%%%%%%%%%%%%%%%%%%%%%%%%%%%%%%%%%%
\begin{theorem}\label{t384696}
Let $2\leq n\leq N$ and 
$\widetilde{p}^{(2)}
=\sum_{R\subseteq\set{n}:\,|R|=2}
\widetilde{p}_{\set{n}\setminus R}
\sum_{j=1}^N\prod_{r\in R}(z_{j,r}-\widetilde{z}_r)$. Then
\begin{align}
\lefteqn{\overline{p}_{\set{n}}
-\frac{N!\,\widetilde{p}_{\set{n}}}{(N-n)!}
+\frac{(N-2)!}{(N-n)!}\widetilde{p}^{(2)}}\nonumber\\
&\qquad=\frac{1}{2N^2}\sum_{k\in\set{n}\setminus\set{2}}h_{k,n}
\sum_{\newatop{R\subseteq\set{n}}{\card{R}=k}}
\sum_{(r,s,t)\in R_{\neq}^3}\widetilde{p}_{\set{n}\setminus R}
\sum_{j\in\set{N}_{\neq}^n}
y_{j_r,j_s,r}y_{j_r,j_s,s}y_{j_r,j_t,t}p_{j,R\setminus\{r,s,t\}}
\nonumber\\
&\qquad\quad{}+\frac{1}{8N^2}\sum_{k\in\set{n}\setminus\set{3}}h_{k,n}
\sum_{\newatop{R\subseteq \set{n}}{\card{R}=k}}
\sum_{(q,r,s,t)\in R_{\neq}^4}\widetilde{p}_{\set{n}\setminus R}
\sum_{j\in\set{N}_{\neq}^n}
y_{j_q,j_r,q}y_{j_q,j_r,r}y_{j_s,j_t,s}y_{j_s,j_t,t}
p_{j,R\setminus\{q,r,s,t\}},\label{e267566}
\end{align}
where 
$h_{k,n}
=\frac{(n+k-2)(n-k+1)}{k(k-1)(k-2)\binomial{n}{k}}$ for 
$k\in\set{n}\setminus\set{2}$. 
If $n=2$ the right-hand side of the equality in 
\eqref{e267566} is defined to be zero. 
\end{theorem}
%%%%%%%%%%%%%%%%%%%%%%%%%%%%%%%%%%%%%%%%%%%%%%%%%%%%%%%%%%%%%%%%%%%%%%
\Proof
We have 
\begin{align}
\widetilde{p}^{(2)}
&=\frac{1}{2}
\sum_{(r,s)\in \set{n}_{\neq}^2}
\sum_{j=1}^N(z_{j,r}-\widetilde{z}_r)(z_{j,s}-\widetilde{z}_s)
\widetilde{p}_{\set{n}\setminus \{r,s\}},\label{e865739}
\end{align}
where, for $(r,s)\in\set{n}_{\neq}^2$, 
\begin{align}
\sum_{j=1}^N(z_{j,r}-\widetilde{z}_r)(z_{j,s}-\widetilde{z}_s)
=\sum_{j\in\set{N}}z_{j,r}z_{j,s}
-N\widetilde{z}_r\widetilde{z}_s
=\frac{1}{2N}\sum_{j\in\set{N}_{\neq}^{2}}y_{j_1,j_2,r}y_{j_1,j_2,s}.
\label{e874628}
\end{align}
Furthermore
\begin{align}\label{e98176}
\frac{1}{2}=\sum_{k=2}^n\frac{k(k-1)}{kn(n-1)}
=\sum_{k=2}^n\frac{\binomial{n-2}{k-2}}{k\binomial{n}{k}}.
\end{align}
In view of \eqref{e865739}, \eqref{e874628} and \eqref{e98176}, 
we see that  
\begin{align}
\frac{(N-2)!}{(N-n)!}\widetilde{p}^{(2)}
&=\frac{(N-2)!}{(N-n)!\,4N}
\sum_{(r,s)\in \set{n}_{\neq}^2}
\sum_{j\in\set{N}_{\neq}^{2}}y_{j_1,j_2,r}y_{j_1,j_2,s} 
\widetilde{p}_{\set{n}\setminus \{r,s\}}\nonumber\\
&=\sum_{k=2}^n\frac{\binomial{n-2}{k-2}}{2Nk\binomial{n}{k}}
\sum_{(r,s)\in \set{n}_{\neq}^2}\sum_{j\in\set{N}_{\neq}^{n}}
y_{j_r,j_{s},r}y_{j_r,j_{s},s} 
\widetilde{p}_{\set{n}\setminus \{r,s\}}\nonumber\\
&=\sum_{k=2}^n\frac{1}{2Nk\binomial{n}{k}}
\sum_{\newatop{R\subseteq\set{n}}{\card{R}=k}}
\sum_{(r,s)\in R_{\neq}^2}\sum_{j\in\set{N}_{\neq}^{n}}
y_{j_r,j_{s},r}y_{j_r,j_{s},s} 
\widetilde{p}_{\set{n}\setminus \{r,s\}}.\label{e7468322}
\end{align}
Combining \eqref{e8366548} and \eqref{e7468322}, we obtain
\begin{align}
\lefteqn{\overline{p}_{\set{n}}
-\frac{N!\,\widetilde{p}_{\set{n}}}{(N-n)!}
+\frac{(N-2)!}{(N-n)!}\widetilde{p}^{(2)}}\nonumber\\
&=-\sum_{k=2}^n\frac{1}{2Nk\binomial{n}{k}}
\sum_{\newatop{R\subseteq\set{n}}{\card{R}=k}}
\sum_{(r,s)\in R_{\neq}^2}\widetilde{p}_{\set{n}\setminus R}
\sum_{j\in\set{N}_{\neq}^{n}}
y_{j_r,j_{s},r}y_{j_r,j_{s},s}(p_{j,R\setminus\{r,s\}}
-\widetilde{p}_{R\setminus \{r,s\}}).\label{e8362795}
\end{align}
In particular, we see that, for $n=2$, \eqref{e267566} is true. 
From now on, let $n\geq 3$. For $k\in\set{n}\setminus\{1\}$, 
$R\subseteq\set{n}$ with $\card{R}=k$,
$(r,s)\in R_{\neq}^2$ and $j\in\set{N}_{\neq}^n$, we have 
\begin{align}
p_{j,R\setminus\{r,s\}}-\widetilde{p}_{R\setminus\{r,s\}}
&=p_{j,R\setminus\{r,s\}}\widetilde{p}_{\emptyset}-p_{j,\emptyset}
\widetilde{p}_{R\setminus\{r,s\}}\nonumber\\
&=\sum_{\ell=1}^{k-2}\Big(
\frac{1}{\binomial{k-2}{\ell}}
\sum_{\newatop{L\subseteq R\setminus\{r,s\}}{\card{L}=\ell}}p_{j,L}
\widetilde{p}_{R\setminus (L\cup\{r,s\})}
-\frac{1}{\binomial{k-2}{\ell-1}}
\sum_{\newatop{L\subseteq R\setminus\{r,s\}}{\card{L}=\ell-1}}
p_{j,L}\widetilde{p}_{R\setminus (L\cup\{r,s\})}\Big)\nonumber\\
&=\sum_{\ell=1}^{k-2}\frac{1}{\binomial{k-2}{\ell}\ell}
\sum_{\newatop{L\subseteq R\setminus\{r,s\}}{\card{L}=\ell}}
\sum_{t\in L}\widetilde{p}_{R\setminus (L\cup\{r,s\})}(p_{j,L}
-\widetilde{z}_{t}p_{j,L\setminus\{t\}}).\label{e487674}
\end{align}
For $R\subseteq\set{n}$ with $\card{R}\geq 3$, 
$(r,s,t)\in R_{\neq}^3$, Lemma \ref{l457638} implies that
\begin{align}
F(r,s,t,R)
&:=\sum_{j\in\set{N}_{\neq}^{n}} y_{j_r,j_{s},r}y_{j_r,j_{s},s}
(p_{j,R\setminus\{r,s\}}-\widetilde{z}_{t}p_{j,R\setminus\{r,s,t\}})
\nonumber\\
&=D_1(r,s,t,R\setminus\{r,s,t\})-D_2(r,s,t,R\setminus\{r,s,t\}),
\label{e137221}
\end{align}
where $D_\nu(r,s,t,R\setminus\{r,s,t\})$ for $\nu\in\set{2}$
are defined as in \eqref{e186573} and \eqref{e186574}. 
Using \eqref{e8362795} and \eqref{e487674}, we get
\begin{align*}
\lefteqn{\overline{p}_{\set{n}}
-\frac{N!\,\widetilde{p}_{\set{n}}}{(N-n)!}
+\frac{(N-2)!}{(N-n)!}\widetilde{p}^{(2)}}\\
&=-\sum_{k=2}^n\sum_{\ell=1}^{k-2}\frac{1}{2Nk\binomial{n}{k}}
\frac{1}{\binomial{k-2}{\ell}\ell}
\sum_{\newatop{R\subseteq\set{n}}{\card{R}=k}}
\sum_{(r,s)\in R_{\neq}^2}
\sum_{\newatop{L\subseteq R\setminus\{r,s\}}{\card{L}=\ell}}
\sum_{t\in L}\widetilde{p}_{\set{n}\setminus (L\cup\{r,s\})}
F(r,s,t,L\cup\{r,s\})\\
&=-\sum_{\ell=1}^{n-2}
\sum_{k=\ell+2}^n\frac{(n-\ell-2)!}{2Nk\ell\binomial{n}{k}
\binomial{k-2}{\ell}(n-k)!(k-\ell-2)!}
\sum_{\newatop{L'\subseteq \set{n}}{\card{L'}=\ell+2}}
\sum_{(r,s,t)\in (L')_{\neq}^3}
\widetilde{p}_{\set{n}\setminus L'}F(r,s,t,L'),
\end{align*}
where, for $\ell\in\set{n-2}$, 
\begin{align*}
\lefteqn{\sum_{k=\ell+2}^n\frac{(n-\ell-2)!}{k\ell\binomial{n}{k}
\binomial{k-2}{\ell}(n-k)!(k-\ell-2)!}
=\frac{(\ell-1)!(n-\ell-2)!}{n!}\sum_{k=\ell+2}^{n}(k-1)}\\
&=\frac{(\ell-1)!(n-\ell-2)!(n+\ell)(n-\ell-1)}{2\,n!} 
=\frac{(n+\ell)(n-\ell-1)}{2\ell(\ell+1)(\ell+2)
\binomial{n}{\ell+2}}=\frac{h_{\ell+2,n}}{2}. 
\end{align*}
Hence 
\begin{align*}
\overline{p}_{\set{n}}
-\frac{N!\,\widetilde{p}_{\set{n}}}{(N-n)!}
+\frac{(N-2)!}{(N-n)!}\widetilde{p}^{(2)}
&=-\frac{1}{4N}\sum_{k=3}^{n}h_{k,n}
\sum_{\newatop{R\subseteq \set{n}}{\card{R}=k}}
\sum_{(r,s,t)\in R_{\neq}^3}
\widetilde{p}_{\set{n}\setminus R}F(r,s,t,R).
\end{align*}
Using this in combination with \eqref{e137221}, the assertion is 
shown. \hfill\qed\medskip

%%%%%%%%%%%%%%%%%%%%%%%%%%%%%%%%%%%%%%%%%%%%%%%%%%%%%%%%%%%%%%%%%%%%%%
As a corollary of Theorem \ref{t384696}, we give a second order 
expansion for the normalized elementary symmetric polynomials. 
%%%%%%%%%%%%%%%%%%%%%%%%%%%%%%%%%%%%%%%%%%%%%%%%%%%%%%%%%%%%%%%%%%%%%%
\begin{corollary}
Let $n,N\in\NN$ with $2\leq n\leq N$, $z_j\in\calz$ for $j\in\set{N}$ 
and $\widetilde{z}=\frac{1}{N}\sum_{j=1}^Nz_j$. 
For $A\subseteq\set{N}$ and $k\in\ZZ$, let
$E_{A,k}=\sum_{\newatop{J\subseteq A}{\card{J}=k}}\prod_{j\in J}z_j$.
Then, we have 
\begin{align} 
\lefteqn{\frac{1}{\binomial{N}{n}}E_{\set{N},n}
-\widetilde{z}^n
+\frac{n(n-1)}{2N(N-1)}
\sum_{j=1}^N(z_{j}-\widetilde{z})^2\widetilde{z}^{n-2}}\nonumber\\
&=\frac{1}{2N^2}\sum_{(r,s,t)\in\set{N}_{\neq}^3}
(z_{r}-z_{s})^2(z_{r}-z_{t})\sum_{k\in\set{n}\setminus\set{2}}
\widetilde{h}_{k,n,N}\widetilde{z}^{n-k}
E_{\set{N}\setminus\{r,s,t\},k-3}\nonumber\\
&\quad{}+\frac{1}{8N^2}\sum_{(q,r,s,t)
\in\set{N}_{\neq}^4}(z_{q}-z_{r})^2(z_{s}-z_{t})^2
\sum_{k\in\set{n}\setminus\set{3}}\widetilde{h}_{k,n,N}
\widetilde{z}^{n-k} E_{\set{N}\setminus\{q,r,s,t\},k-4},
\label{e365765}
\end{align}
where 
$\widetilde{h}_{k,n,N}
=\frac{(n+k-2)(n-k+1)}{k(k-1)(k-2)\binomial{N}{k}}$
for $k\in\set{n}\setminus\set{2}$. 
If $n=2$, the right-hand side of the equality in 
\eqref{e365765} is defined to be zero. 
\end{corollary}
%%%%%%%%%%%%%%%%%%%%%%%%%%%%%%%%%%%%%%%%%%%%%%%%%%%%%%%%%%%%%%%%%%%%%%
\Proof
Similarly as in the proof of Corollary \ref{c526745}, 
Identity \eqref{e365765} follows from 
Theorem \ref{t384696} in the case 
that $Z$ has identical columns. Indeed, letting 
$z_{j,1}=\dots=z_{j,n}=z_j$ for all $j\in\set{N}$, then
$\widetilde{z}_1=\dots=\widetilde{z}_n=\widetilde{z}$ and
\begin{align*}
\overline{p}_{\set{n}}
=n!E_{\set{N},n}, 
\qquad
\widetilde{p}_{\set{n}}
=\widetilde{z}^n,\qquad 
\widetilde{p}^{(2)}
=\binomial{n}{2}
\sum_{j=1}^N(z_{j}-\widetilde{z})^2\widetilde{z}^{n-2}.
\end{align*}
Therefore 
\begin{align}
\frac{1}{\binomial{N}{n}}E_{\set{N},n}-\widetilde{z}^n
+\frac{n(n-1)}{2N(N-1)}
\sum_{j=1}^N(z_{j}-\widetilde{z})^2\widetilde{z}^{n-2}
&=\frac{(N-n)!}{N!}
\Big(\overline{p}_{\set{n}}
-\frac{N!\,\widetilde{p}_{\set{n}}}{(N-n)!}
+\frac{(N-2)!}{(N-n)!}\widetilde{p}^{(2)}\Big)
\nonumber\\
&=\frac{(N-n)!}{N!}\Bigl(\frac{M_1}{2N^2} +\frac{M_2}{8N^2}\Bigr), 
\label{e578581}
\end{align} 
where 
\begin{align*}
M_1&=\sum_{k\in\set{n}\setminus\set{2}}h_{k,n}
\sum_{\newatop{R\subseteq\set{n}}{\card{R}=k}}
\sum_{(r,s,t)\in R_{\neq}^3}\widetilde{p}_{\set{n}\setminus R}
\sum_{j\in\set{N}_{\neq}^n}
y_{j_r,j_s,r}y_{j_r,j_s,s}y_{j_r,j_t,t}p_{j,R\setminus\{r,s,t\}},\\
M_2&=\sum_{k\in\set{n}\setminus\set{3}}h_{k,n}
\sum_{\newatop{R\subseteq \set{n}}{\card{R}=k}}
\sum_{(q,r,s,t)\in R_{\neq}^4}\widetilde{p}_{\set{n}\setminus R}
\sum_{j\in\set{N}_{\neq}^n}
y_{j_q,j_r,q}y_{j_q,j_r,r}y_{j_s,j_t,s}y_{j_s,j_t,t}
p_{j,R\setminus\{q,r,s,t\}}.
\end{align*}
Here 
\begin{align*}
M_1
&=\sum_{k\in\set{n}\setminus\set{2}}h_{k,n}
\sum_{\newatop{R\subseteq\set{n}}{\card{R}=k}}
\sum_{(r,s,t)\in R_{\neq}^3}\widetilde{z}^{n-k}
\sum_{j\in\set{N}_{\neq}^n}
(z_{j_r}-z_{j_s})^2
(z_{j_r}-z_{j_t})\prod_{\ell\in R\setminus\{r,s,t\}}z_{j_\ell}\\
&=\sum_{(u,v,w)\in\set{N}_{\neq}^3}
(z_{u}-z_{v})^2(z_{u}-z_{w})\sum_{k\in\set{n}\setminus\set{2}}
h_{k,n}\widetilde{z}^{n-k}
\sum_{\newatop{R\subseteq\set{n}}{\card{R}=k}}
\sum_{(r,s,t)\in R_{\neq}^3}
\sum_{\newatop{j\in\set{N}_{\neq}^n}{j_r=u,j_s=v,j_t=w}}
\prod_{\ell\in R\setminus\{r,s,t\}}z_{j_\ell},
\end{align*}
where, for $(u,v,w)\in\set{N}_{\neq}^3$ and 
$k\in\set{n}\setminus\set{2}$, 
\begin{align*}
\sum_{\newatop{R\subseteq\set{n}}{\card{R}=k}}
\sum_{(r,s,t)\in R_{\neq}^3}
\sum_{\newatop{j\in\set{N}_{\neq}^n}{j_r=u,j_s=v,j_t=w}}
\prod_{\ell\in R\setminus\{r,s,t\}}z_{j_\ell}
=\frac{N!\binomial{n}{k}}{(N-n)!\binomial{N}{k}} 
E_{\set{N}\setminus\{u,v,w\},k-3}
\end{align*}
and 
$h_{k,n}\frac{\binomial{n}{k}}{\binomial{N}{k}}
=\widetilde{h}_{k,n,N}$. 
This implies that 
\begin{align} % number of formula is necessary !
M_1
&=\frac{N!}{(N-n)!}\sum_{(r,s,t)\in\set{N}_{\neq}^3}
(z_{r}-z_{s})^2(z_{r}-z_{t})\sum_{k\in\set{n}\setminus\set{2}}
\widetilde{h}_{k,n,N}\widetilde{z}^{n-k}
E_{\set{N}\setminus\{r,s,t\},k-3}. 
\end{align}
Furthermore
\begin{align*}
M_2
&=\sum_{k\in\set{n}\setminus\set{3}}h_{k,n}
\sum_{\newatop{R\subseteq \set{n}}{\card{R}=k}}
\sum_{(q,r,s,t)\in R_{\neq}^4}\widetilde{z}^{n-k}
\sum_{j\in\set{N}_{\neq}^n}
(z_{j_q}-z_{j_r})^2(z_{j_s}-z_{j_t})^2
\prod_{\ell\in R\setminus\{q,r,s,t\}} z_{j_\ell}\\
&=\sum_{(u,v,w,x)\in\set{N}_{\neq}^4}(z_{u}-z_{v})^2(z_{w}-z_{x})^2\\
&\quad{}\times\sum_{k\in\set{n}\setminus\set{3}}
h_{k,n}\widetilde{z}^{n-k}
\sum_{\newatop{R\subseteq \set{n}}{\card{R}=k}}
\sum_{(q,r,s,t)\in R_{\neq}^4}
\sum_{\newatop{j\in\set{N}_{\neq}^n}{j_q=u,j_r=v,j_s=w,j_t=x}}
\prod_{\ell\in R\setminus\{q,r,s,t\}} z_{j_\ell},
\end{align*}
where, for $(u,v,w,x)\in\set{N}_{\neq}^4$ and  
$k\in\set{n}\setminus\set{3}$,
\begin{align*}
\sum_{\newatop{R\subseteq \set{n}}{\card{R}=k}}
\sum_{(q,r,s,t)\in R_{\neq}^4}
\sum_{\newatop{j\in\set{N}_{\neq}^n}{j_q=u,j_r=v,j_s=w,j_t=x}}
\prod_{\ell\in R\setminus\{q,r,s,t\}} z_{j_\ell}
=\frac{N!\binomial{n}{k}}{(N-n)!\binomial{N}{k}}
E_{\set{N}\setminus\{u,v,w,x\},k-4}.
\end{align*}
Hence
\begin{align}
M_2
&=\frac{N!}{(N-n)!}\sum_{(q,r,s,t)
\in\set{N}_{\neq}^4}(z_{q}-z_{r})^2(z_{s}-z_{t})^2
\sum_{k\in\set{n}\setminus\set{3}}
\widetilde{h}_{k,n,N}\widetilde{z}^{n-k}
E_{\set{N}\setminus\{q,r,s,t\},k-4}.\label{e578583}
\end{align}
Combining \eqref{e578581}--\eqref{e578583}, the assertion is 
shown. \hfill\qed
%%%%%%%%%%%%%%%%%%%%%%%%%%%%%%%%%%%%%%%%%%%%%%%%%%%%%%%%%%%%%%%%%%%%%%

\section{Approximation of normalized permanents}\label{s749687}
%%%%%%%%%%%%%%%%%%%%%%%%%%%%%%%%%%%%%%%%%%%%%%%%%%%%%%%%%%%%%%%%%%%%%%
In this section, we employ the notation of Section \ref{s326765} 
with $\calz:=\CC$. 
It should be mentioned that, unless stated otherwise, we do not 
assume that the numbers $|z_{j,r}|$ for $j\in\set{N}$, $r\in\set{n}$ 
are bounded by one. 
%%%%%%%%%%%%%%%%%%%%%%%%%%%%%%%%%%%%%%%%%%%%%%%%%%%%%%%%%%%%%%%%%%%%%%
\subsection{Main approximation results}\label{s4518u5}
%%%%%%%%%%%%%%%%%%%%%%%%%%%%%%%%%%%%%%%%%%%%%%%%%%%%%%%%%%%%%%%%%%%%%%
The first results in this section are Theorems \ref{t74739} 
and \ref{t4216459} below, the proof of which require the following 
lemma. 
%%%%%%%%%%%%%%%%%%%%%%%%%%%%%%%%%%%%%%%%%%%%%%%%%%%%%%%%%%%%%%%%%%%%%%
\begin{lemma} \label{l378687}
\begin{enumerate}

\item\label{l378687.a}%
(Hadamard type permanent inequality)
Without any further restrictions, we have 
\begin{align*}
|\Per(Z)|
\leq \frac{N!}{(N-n)!}\prod_{r=1}^n\Big(
\frac{1}{N}\sum_{j=1}^N|z_{j,r}|^2\Big)^{1/2}.
\end{align*}

\item\label{l378687.b}%
(Br\'{e}gman-Minc permanent inequality)
If $Z\in\{0,1\}^{N\times n}$, then 
\begin{align*}
\Per(Z)\leq \frac{N!}{(N-n)!}
\prod_{r=1}^n
\frac{\zeta(N\widetilde{z}_r)}{(N!)^{1/N}},
\end{align*}
where $\zeta(k)=(k!)^{1/k}$ for $k\in\NN$ and $\zeta(0)=0$. 
\end{enumerate}

\end{lemma} 
%%%%%%%%%%%%%%%%%%%%%%%%%%%%%%%%%%%%%%%%%%%%%%%%%%%%%%%%%%%%%%%%%%%%%%
\Proof Let us first consider the case $n=N$. For \eqref{l378687.a}, 
see Carlen et al.\ \citet[Theorem 1.1]{MR2275869} or 
Cobos et al.\ \citet[Theorem 5.1]{MR2256997}.
As stated in \cite[Introduction]{MR2275869}, this can also be 
obtained from 
Theorem~9.1.1 in Appendix~1 of Nesterov and Nemirovskii 
\citet{MR1258086}. 
Part \eqref{l378687.b} was conjectured by Minc \citet{MR0155843} and 
proved by Br\`{e}gman \citet{MR0327788}. 
The general case of $n\in\set{N}$ follows from the above and the 
simple observation that $\Per(Z)=\frac{\Per(Z')}{(N-n)!}$, where 
$Z'=(z_{j,r}')\in\CC^{N\times N}$ with $z'_{j,r}=z_{j,r}$ for 
$j\in\set{N}$, $r\in\set{n}$ and $z'_{j,r}=1$ for $j\in\set{N}$, 
$r\in\set{N}\setminus\set{n}$. \hfill\qed
\medskip

%%%%%%%%%%%%%%%%%%%%%%%%%%%%%%%%%%%%%%%%%%%%%%%%%%%%%%%%%%%%%%%%%%%%%%
We note that in \cite[Lemma 2.2]{MR3314090}, a Hadamard type 
inequality for the permanent of a matrix with zero column sums was 
shown, which is uniformly better than the general bound in 
Lemma \ref{l378687}\eqref{l378687.a}. 

The next theorem contains an improvement of the 
inequalities \eqref{e6285}--\eqref{e63286508} from the introduction. 

%%%%%%%%%%%%%%%%%%%%%%%%%%%%%%%%%%%%%%%%%%%%%%%%%%%%%%%%%%%%%%%%%%%%%%
\begin{theorem}\label{t74739}
Let us assume that $2\leq n\leq N$ and set
\begin{gather*}
\beta=\frac{1}{n}\sum_{r=1}^n|\widetilde{z}_r|^2,\qquad
\vartheta=\frac{1}{N(N-1)\sqrt{n(n-1)}}
\Big(\sum_{(r,s)\in \set{n}_{\neq}^2}
\Big(\sum_{(u,v)\in\set{N}_{\neq}^2}|y_{u,v,r}y_{u,v,s}|
\Big)^2\Big)^{1/2},\\
\kappa
= 
\begin{cases}
\displaystyle
\frac{1}{(n-2)(N-2)}
\max_{(u,v)\in\set{N}_{\neq}^2,(r,s)\in\set{n}_{\neq}^2}
\sum_{j\in\set{N}\setminus\{u,v\}}
\sum_{\ell\in\set{n}\setminus\{r,s\}}|z_{j,\ell}|^2, &
\quad\mbox{if } n\geq 3,\\
1, &\quad\mbox{if } n=2,
\end{cases} \\
f_n(x_1,x_2)=\sum_{k=2}^n(k-1)x_1^{n-k}x_2^{k-2}
=\frac{(n-1)x_2^{n}-nx_1x_2^{n-1}+x_1^n}{(x_2-x_1)^2},
\qquad(x_1,x_2\in\RR).
\end{gather*}
Then
\begin{align} \label{e7465283}
\Big|\frac{(N-n)!}{N!}\overline{p}_{\set{n}}
-\widetilde{p}_{\set{n}}\Big|
&\leq\frac{\vartheta }{2N}f_n(\sqrt{\beta},\sqrt{\kappa}).
\end{align}
\end{theorem}
%%%%%%%%%%%%%%%%%%%%%%%%%%%%%%%%%%%%%%%%%%%%%%%%%%%%%%%%%%%%%%%%%%%%%%
\Proof 
In the case $n=2$, \eqref{e5703} gives
\begin{align}
\frac{(N-n)!}{N!}\overline{p}_{\set{n}}-\widetilde{p}_{\set{n}}
=-\frac{1}{2N^2(N-1)}\sum_{(u,v)\in\set{N}_{\neq}^2}
y_{u,v,1}y_{u,v,2},\label{e51764386}
\end{align}
which together with the identities
$\vartheta
=\frac{1}{N(N-1) }
\sum_{(u,v)\in\set{N}_{\neq}^2}|y_{u,v,1}y_{u,v,2}|$ and
$f_n(\sqrt{\beta},\sqrt{\kappa})=1$
implies \eqref{e7465283}. Let us now assume that $3\leq n\leq N$. 
From \eqref{e8366548}, we get
\begin{align}
\lefteqn{\Big|\frac{(N-n)!}{N!}\overline{p}_{\set{n}}
-\widetilde{p}_{\set{n}}\Big|
\leq \frac{1}{N!}\sum_{k=2}^n\frac{(N-n)!}{2Nk\binomial{n}{k}}
\sum_{\newatop{R\subseteq\set{n}}{\card{R}=k}}
\sum_{(r,s)\in R_{\neq}^2}\sum_{j\in\set{N}_{\neq}^{n}}
|y_{j_r,j_{s},r}y_{j_r,j_{s},s}p_{j,R\setminus\{r,s\}}
\widetilde{p}_{\set{n}\setminus R}|}\nonumber\\
&=\frac{1}{N!}\sum_{k=2}^n\frac{(N-k)!}{2Nk\binomial{n}{k}}
\sum_{\newatop{R\subseteq\set{n}}{\card{R}=k}}
\sum_{(r,s)\in R_{\neq}^2}|\widetilde{p}_{\set{n}\setminus R}|
\sum_{(u,v)\in\set{N}_{\neq}^2}|y_{u,v,r}y_{u,v,s}|
\sum_{j\in(\set{N}\setminus\{u,v\})_{\neq}^{R\setminus\{r,s\}}}
|p_{j,R\setminus\{r,s\}}|.\label{e615496}
\end{align}
Lemma \ref{l378687}\eqref{l378687.a} implies that, for 
$R\subseteq\set{n}$ with $\card{R}=k\geq 2$, $(r,s)\in R_{\neq}^2$, 
$(u,v)\in\set{N}_{\neq}^2$,
\begin{align*}
\sum_{j\in(\set{N}\setminus\{u,v\})_{\neq}^{R\setminus\{r,s\}}}
|p_{j,R\setminus\{r,s\}}|
&\leq \frac{(N-2)!}{(N-k)!}
\prod_{\ell\in R\setminus\{r,s\}}\kappa_{u,v,\ell}^{1/2},
\end{align*}
where  $\kappa_{u,v,\ell}
=\frac{1}{N-2}\sum_{j\in\set{N}\setminus\{u,v\}}|z_{j,\ell}|^2$
for $(u,v)\in\set{N}_{\neq}^2$ and $\ell\in\set{n}$.
For $k\in\set{n}\setminus\set{1}$, the 
Cauchy-Schwarz inequality gives
\begin{align*}
\lefteqn{\sum_{\newatop{R\subseteq\set{n}}{\card{R}=k}}
\sum_{(r,s)\in R_{\neq}^2}|\widetilde{p}_{\set{n}\setminus R}|
\sum_{(u,v)\in\set{N}_{\neq}^2}|y_{u,v,r}y_{u,v,s}|
\prod_{\ell\in R\setminus\{r,s\}}\kappa_{u,v,\ell}^{1/2}}\\
&\leq\Big(\sum_{\newatop{R\subseteq\set{n}}{\card{R}=k}}
\sum_{(r,s)\in R_{\neq}^2}|\widetilde{p}_{\set{n}\setminus R}|^2
\Big)^{1/2}
\Big(\sum_{(r,s)\in \set{n}_{\neq}^2}
\sum_{\newatop{R\subseteq\set{n}\setminus\{r,s\}}{\card{R}=k-2}}
\Big(\sum_{(u,v)\in\set{N}_{\neq}^2}|y_{u,v,r}y_{u,v,s}|
\prod_{\ell\in R}\kappa_{u,v,\ell}^{1/2}\Big)^2\Big)^{1/2}.
\end{align*}
Using Maclaurin's inequality 
(see Hardy et al.\ \citet[Theorem 52, page 52]{MR0046395}), we obtain 
\begin{align*}
\sum_{\newatop{R\subseteq\set{n}}{\card{R}=k}}
\sum_{(r,s)\in R_{\neq}^2}|\widetilde{p}_{\set{n}\setminus R}|^2
=k(k-1)\sum_{\newatop{R\subseteq\set{n}}{\card{R}=n-k}}
|\widetilde{p}_{R}|^2
\leq k(k-1)\binomial{n}{n-k}\beta^{n-k}.
\end{align*}
Furthermore, for $(r,s)\in\set{n}_{\neq}^2$,  
\begin{align*}
\lefteqn{\sum_{\newatop{R\subseteq\set{n}\setminus\{r,s\}}{
\card{R}=k-2}}
\Big(\sum_{(u,v)\in\set{N}_{\neq}^2}|y_{u,v,r}y_{u,v,s}|
\prod_{\ell\in R}\kappa_{u,v,\ell}^{1/2}\Big)^2}\\
&=\sum_{(u,v)\in\set{N}_{\neq}^2}
\sum_{(u',v')\in\set{N}_{\neq}^2}|y_{u,v,r}y_{u,v,s}
y_{u',v',r}y_{u',v',s}|
\sum_{\newatop{R\subseteq\set{n}\setminus\{r,s\}}{\card{R}=k-2}}
\prod_{\ell\in R}(\kappa_{u,v,\ell}\kappa_{u',v',\ell})^{1/2},
\end{align*}
where, for $(u,v),(u',v')\in\set{N}_{\neq}^2$, 
the Cauchy-Schwarz inequality gives 
\begin{align*}
\sum_{\newatop{R\subseteq\set{n}\setminus\{r,s\}}{\card{R}=k-2}}
\prod_{\ell\in R}(\kappa_{u,v,\ell}\kappa_{u',v',\ell})^{1/2}
\leq \Big(\sum_{\newatop{R\subseteq\set{n}\setminus\{r,s\}}{
\card{R}=k-2}}
\prod_{\ell\in R}\kappa_{u,v,\ell}\Big)^{1/2}
\Big(\sum_{\newatop{R\subseteq\set{n}\setminus\{r,s\}}{\card{R}=k-2}}
\prod_{\ell\in R}\kappa_{u',v',\ell}\Big)^{1/2}
\end{align*}
and Maclaurin's inequality implies that 
\begin{gather*}
\sum_{\newatop{R\subseteq\set{n}\setminus\{r,s\}}{
\card{R}=k-2}}\prod_{\ell\in R}\kappa_{u,v,\ell}
\leq \binomial{n-2}{k-2}\Big(\frac{1}{n-2}
\sum_{\ell\in\set{n}\setminus\{r,s\}}\kappa_{u,v,\ell}\Big)^{k-2}
\leq \binomial{n-2}{k-2}\kappa^{k-2}.
\end{gather*}
Hence
\begin{align*}
\sum_{(r,s)\in\set{n}_{\neq}^2}
\sum_{\newatop{R\subseteq\set{n}\setminus\{r,s\}}{\card{R}=k-2}}
\Big(\sum_{(u,v)\in\set{N}_{\neq}^2}|y_{u,v,r}y_{u,v,s}|
\prod_{\ell\in R}\kappa_{u,v,\ell}^{1/2}\Big)^2
&\leq \binomial{n-2}{k-2}\kappa^{k-2}
n(n-1)(N(N-1)\vartheta)^2.
\end{align*}
Combining the inequalities above, we obtain 
\begin{align*}
\Big|\frac{(N-n)!}{N!}\overline{p}_{\set{n}}
-\widetilde{p}_{\set{n}}\Big|
&\leq \frac{1}{N!}\sum_{k=2}^n\frac{(N-2)!}{2Nk\binomial{n}{k}}
\Big(\sum_{\newatop{R\subseteq\set{n}}{\card{R}=k}}
\sum_{(r,s)\in R_{\neq}^2}|\widetilde{p}_{\set{n}\setminus R}|^2
\Big)^{1/2}\\
&\quad{}\times \Big(\sum_{(r,s)\in \set{n}_{\neq}^2}
\sum_{\newatop{R\subseteq\set{n}\setminus\{r,s\}}{\card{R}=k-2}}
\Big(\sum_{(u,v)\in\set{N}_{\neq}^2}|y_{u,v,r}y_{u,v,s}|
\prod_{\ell\in R}\kappa_{u,v,\ell}^{1/2}\Big)^2\Big)^{1/2}\\
&\leq \frac{\vartheta}{2N}
\sum_{k=2}^n(k-1)\beta^{(n-k)/2}\kappa^{(k-2)/2}\\
&=\frac{\vartheta }{2N}f_n(\sqrt{\beta},\sqrt{\kappa}),
\end{align*}
which implies the assertion. \hfill\qed
\medskip

%%%%%%%%%%%%%%%%%%%%%%%%%%%%%%%%%%%%%%%%%%%%%%%%%%%%%%%%%%%%%%%%%%%%%%
The proof of Theorem \ref{t74739} requires only Part 
\eqref{l378687.a} of Lemma \ref{l378687}. But if $Z$ is a 
$(0,1)$-matrix, then Part \eqref{l378687.b} can also be applied, 
as is shown in the following theorem. 

%%%%%%%%%%%%%%%%%%%%%%%%%%%%%%%%%%%%%%%%%%%%%%%%%%%%%%%%%%%%%%%%%%%%%%
\begin{theorem}\label{t4216459}
Let the assumptions of Theorem \ref{t74739} be valid, where 
we assume that $Z\in\{0,1\}^{N\times n}$. Further, 
let $\zeta$ be as in  Lemma \ref{l378687}\eqref{l378687.b}, 
$\eta_{u,v,\ell}
=\sum_{j\in\set{N}\setminus\{u,v\}}z_{j,\ell}$,
$\widetilde{\kappa}_{u,v,\ell}
=\frac{(\zeta(\eta_{u,v,\ell}))^2}{((N-2)!)^{2/(N-2)}}$
for $(u,v)\in\set{N}_{\neq}^2$, $\ell\in\set{n}$ and set 
\begin{align*}
\widetilde{\kappa}
&=
\begin{cases}
\displaystyle\frac{1}{n-2}
\max_{(u,v)\in\set{N}_{\neq}^2,(r,s)\in\set{n}_{\neq}^2}
\sum_{\ell\in\set{n}\setminus\{r,s\}}
\widetilde{\kappa}_{u,v,\ell},&
\quad\mbox{ if } n\geq 3,\\
1, &\quad\mbox{ if } n= 2. 
\end{cases}
\end{align*}
Then 
\begin{align} \label{e7465284}
\Big|\frac{(N-n)!}{N!}\overline{p}_{\set{n}}
-\widetilde{p}_{\set{n}}\Big|
&\leq\frac{\vartheta }{2N}
f_n(\sqrt{\beta},\sqrt{\widetilde{\kappa}}).
\end{align}
\end{theorem}
%%%%%%%%%%%%%%%%%%%%%%%%%%%%%%%%%%%%%%%%%%%%%%%%%%%%%%%%%%%%%%%%%%%%%%
\Proof 
If $n=2$, the assertion directly follows from Theorem \ref{t74739}. 
In the case $3\leq n\leq N$, the proof is very similar to the one of 
Theorem \ref{t74739}. Indeed, we use \eqref{e615496} 
in combination with Lemma \ref{l378687}\eqref{l378687.b}, which  
implies that, for $R\subseteq\set{n}$ with 
$\card{R}=k\geq 2$, $(r,s)\in R_{\neq}^2$, $(u,v)\in\set{N}_{\neq}^2$,
\begin{align*}
\sum_{j\in(\set{N}\setminus\{u,v\})_{\neq}^{R\setminus\{r,s\}}}
|p_{j,R\setminus\{r,s\}}|
&\leq\frac{(N-2)!}{(N-k)!}\prod_{\ell\in R\setminus\{r,s\}}
\widetilde{\kappa}_{u,v,\ell}^{1/2}, 
\end{align*}
where 
$\frac{1}{n-2}\sum_{\ell\in\set{n}\setminus\{r,s\}}
\widetilde{\kappa}_{u,v,\ell}
\leq \widetilde{\kappa}$.\hfill\qed  \medskip

%%%%%%%%%%%%%%%%%%%%%%%%%%%%%%%%%%%%%%%%%%%%%%%%%%%%%%%%%%%%%%%%%%%%%%

The right-hand sides of \eqref{e7465283} and \eqref{e7465284} 
can be further estimated by using the following lemma. 
The inequalities given there can be used for the 
comparison with the bounds given in the introduction. 
%%%%%%%%%%%%%%%%%%%%%%%%%%%%%%%%%%%%%%%%%%%%%%%%%%%%%%%%%%%%%%%%%%%%%%

\begin{lemma}\label{l1436530}
Let us assume that $2\leq n\leq N$. Let 
$\beta$, $\vartheta$, $f_n$ be as in Theorem \ref{t74739}
and set 
$\alpha
=\frac{1}{nN}\sum_{r\in\set{n}}\sum_{j\in\set{N}}
|z_{j,r}-\widetilde{z}_r|^2$. 
Then 
\begin{gather}
\alpha=\frac{1}{nN}\sum_{j=1}^N\sum_{r=1}^n|z_{j,r}|^2-\beta,
\label{e624371}\\
\vartheta  
\leq \frac{2N\alpha}{N-1},\label{e624372}\\
f_n(x,1)
\leq (n-1)\frac{1-x^{n/2}}{1-x}
\leq (n-1)\min\Big\{\frac{n}{2},\,\frac{1}{1-x}\Big\}, 
\quad (x\in[0,1]).\label{e624373}
\end{gather}

\end{lemma}
%%%%%%%%%%%%%%%%%%%%%%%%%%%%%%%%%%%%%%%%%%%%%%%%%%%%%%%%%%%%%%%%%%%%%%
\Proof
As already mentioned in \cite[Remark 2.9]{MR3314090}, 
\eqref{e624371} is true. 
A repeated application of the Cauchy-Schwarz inequality yields
\begin{align*}
\alpha
&=\frac{1}{2nN^2}
\sum_{r\in \set{n}}\sum_{(u,v)\in \set{N}_{\neq}^2}|y_{u,v,r}|^2
\leq\frac{1}{2N^2\sqrt{n}}
\Big(\sum_{r\in\set{n}}\Big(\sum_{(u,v)\in\set{N}_{\neq}^2}
|y_{u,v,r}|^2\Big)^2\Big)^{1/2}
\end{align*}
and 
\begin{align*}
N^2 (N-1)^2n(n-1)\vartheta^2
&\leq
\sum_{(r,s)\in\set{n}_{\neq}^2}
\Big(\sum_{(u,v)\in\set{N}_{\neq}^2}|y_{u,v,r}|^2\Big)
\sum_{(u,v)\in\set{N}_{\neq}^2}|y_{u,v,s}|^2\\
&=\Big(\sum_{r\in\set{n}}
\sum_{(u,v)\in\set{N}_{\neq}^2}|y_{u,v,r}|^2\Big)^2
-\sum_{r\in\set{n}}\Big(\sum_{(u,v)\in\set{N}_{\neq}^2}
|y_{u,v,r}|^2\Big)^2\\
&\leq 4N^4n^2\alpha^2-4N^4n\alpha^2
=4N^4n(n-1)\alpha^2,
\end{align*}
giving \eqref{e624372}. 
Finally, \eqref{e624373} follows from the Jensen inequality. Indeed, 
since $x\in[0,1]$,
\begin{align*}
f_n(x,1)=\frac{n-1-nx+x^n}{(1-x)^2}
&=\frac{n-1}{1-x}\Big(1-\frac{1}{n-1}\sum_{m=1}^{n-1}x^m\Big)
\leq (n-1)\frac{1-x^{n/2}}{1-x},
\end{align*}
where, for $n\geq 2$, 
$\frac{1-x^{n/2}}{1-x}$ is increasing in $x\in[0,1]$, so that 
$\frac{1-x^{n/2}}{1-x}
=\sum_{m=0}^{n-1}\frac{x^{m/2}}{1+\sqrt{x}}
\leq \frac{n}{2}$. \hfill\qed 

%%%%%%%%%%%%%%%%%%%%%%%%%%%%%%%%%%%%%%%%%%%%%%%%%%%%%%%%%%%%%%%%%%%%%%
\begin{corollary} 
If $2\leq n\leq N$, $|z_{j,r}|\leq1$ for all $j\in\set{N}$, 
$r\in\set{n}$ and $\alpha$, $\beta$, $\vartheta$ are as in 
Theorem \ref{t74739} and Lemma \ref{l1436530}, respectively, then
\begin{align}
\Big|\frac{(N-n)!}{N!}\overline{p}_{\set{n}}
-\widetilde{p}_{\set{n}}\Big|
&\leq \frac{n-1}{2N}\vartheta
\frac{1-\beta^{n/4}}{1-\sqrt{\beta}}\label{e739065}\\
&\leq \frac{n-1}{N-1}\alpha
\min\Big\{\frac{n}{2},\,\frac{1}{1-\sqrt{\beta}}\Big\}
\label{e514385}\\
&\leq (1+\sqrt{\beta})\frac{n-1}{N-1}. \label{e627867}
\end{align}
\end{corollary}
%%%%%%%%%%%%%%%%%%%%%%%%%%%%%%%%%%%%%%%%%%%%%%%%%%%%%%%%%%%%%%%%%%%%%%
\Proof
This can easily be derived from Theorem \ref{t74739}, 
Lemma \ref{l1436530}, and the trivial fact that $\kappa\leq 1$, i.e.\
$f_n(\sqrt{\beta},\sqrt{\kappa})\leq f_n(\sqrt{\beta},1)$. 
\hfill\qed\medskip
%%%%%%%%%%%%%%%%%%%%%%%%%%%%%%%%%%%%%%%%%%%%%%%%%%%%%%%%%%%%%%%%%%%%%%

\begin{remark}\label{r732860}
Inequality \eqref{e514385} implies that, in \eqref{e62865098}, 
the right-hand side can be replaced with the expression 
$(1+\sqrt{\beta})\gamma(\frac{1}{2})$. 
\end{remark}
%%%%%%%%%%%%%%%%%%%%%%%%%%%%%%%%%%%%%%%%%%%%%%%%%%%%%%%%%%%%%%%%%%%%%%

We now discuss the sharpness of some of the inequalities above. 
%%%%%%%%%%%%%%%%%%%%%%%%%%%%%%%%%%%%%%%%%%%%%%%%%%%%%%%%%%%%%%%%%%%%%%

\begin{remark} 
Let us assume that $2=n\leq N$ and that 
$Z=(z_{j,r})\in\RR^{N\times 2}$. 
Then, in \eqref{e624373}, equality holds. Below, we will 
additionally assume the validity of some of the following conditions: 
\begin{gather}
z_{j,1}=z_{j,2} \mbox{ for all } j\in\set{N},\label{e632986}\\
Z \mbox{ has decreasing columns}.
\label{e632987}
\end{gather}
\begin{enumerate}

\item 
If \eqref{e632986} is satisfied, then, in \eqref{e624372}, equality 
holds. 

\item  
The right-hand side of \eqref{e7465283} is equal to 
$\frac{\vartheta}{2N}=\frac{1}{2N^2(N-1)}
\sum_{(u,v)\in\set{N}_{\neq}^2}|y_{u,v,1}y_{u,v,2}|$.
In view of \eqref{e51764386}, we see that, if 
one of the conditions \eqref{e632986} or \eqref{e632987} is satisfied,
then $y_{u,v,1}y_{u,v,2}\geq 0$ such that, 
in \eqref{e7465283}, equality holds. 

\item Let us assume that $Z\in[-1,1]^{N\times 2}$.
If one of the conditions \eqref{e632986} or \eqref{e632987} is 
satisfied, then, in \eqref{e739065}, equality holds. 
If \eqref{e632986} holds, 
then, in \eqref{e514385}, equality holds. 
If additionally $\alpha=1$ and $\beta=0$, then, 
in \eqref{e627867}, equality holds. 

\end{enumerate}

\end{remark}

%%%%%%%%%%%%%%%%%%%%%%%%%%%%%%%%%%%%%%%%%%%%%%%%%%%%%%%%%%%%%%%%%%%%%%
\begin{corollary}\label{c6386}
Let $n,N\in\NN$ with $2\leq n\leq N$, 
$z_j\in\CC$ with $|z_j|\leq 1$ for all $j\in\set{N}$. 
Set $\widetilde{z}=\frac{1}{N}\sum_{j=1}^Nz_j$,
$\kappa
=\frac{1}{N-2}\max_{(u,v)\in\set{N}_{\neq}^2}
\sum_{j\in\set{N}\setminus\{u,v\}}|z_{j}|^2$ if $n\geq 3$,
and $\kappa=1$ otherwise. Let $E_{\set{N},n}
=\sum_{\newatop{J\subseteq \set{N}}{\card{J}=n}}\prod_{j\in J}z_j$.
Then
\begin{align}
\Big|\frac{1}{\binomial{N}{n}}
E_{\set{N},n}
-\widetilde{z}^n\Big|
&\leq \frac{f_n(|\widetilde{z}|,\sqrt{\kappa})}{N(N-1)}
\sum_{j=1}^N|z_j-\widetilde{z}|^2\label{e6521765}\\
&\leq \frac{n(n-1)}{N(N-1)}\sum_{j=1}^N|z_j-\widetilde{z}|^2
\min\Big\{\frac{1}{2},\,\frac{1}{n(1-|\widetilde{z}|)}\Big\}.
\label{e6521766}
\end{align}
\end{corollary}
%%%%%%%%%%%%%%%%%%%%%%%%%%%%%%%%%%%%%%%%%%%%%%%%%%%%%%%%%%%%%%%%%%%%%%
\Proof 
Let us consider the matrix $Z$ with 
$z_{j,r}=z_j$ for all $j\in\set{N}$ and $r\in\set{n}$, i.e.\ 
$Z$ has identical columns. Using the notation in 
Theorem~\ref{t74739}, we obtain
\begin{gather*}
\frac{1}{\binomial{N}{n}}E_{\set{N},n}
=\frac{(N-n)!}{N!}\overline{p}_{\set{n}}, \qquad 
\widetilde{z}^n=\widetilde{p}_{\set{n}},\qquad
\beta=|\widetilde{z}|^2,\\
\frac{(N-1) \vartheta}{2}
=\frac{1}{2N}\sum_{(u,v)\in\set{N}_{\neq}^2}|z_{u}-z_v|^2
=\sum_{j=1}^N|z_{j}-\widetilde{z}|^2
\end{gather*}
and hence 
\begin{align*}
\Big|\frac{1}{\binomial{N}{n}}E_{\set{N},n}
-\widetilde{z}^n\Big|
&=\Big|\frac{(N-n)!}{N!}\overline{p}_{\set{n}}
-\widetilde{p}_{\set{n}}\Big|
\leq \frac{\vartheta }{2N}f_n(|\widetilde{z}|,\sqrt{\kappa})
=\frac{f_n(|\widetilde{z}|,\sqrt{\kappa})}{N(N-1)}
\sum_{j=1}^N|z_{j}-\widetilde{z}|^2,
\end{align*}
which proves \eqref{e6521765}. 
Inequality \eqref{e6521766} follows with the help of \eqref{e624373}. 
\hfill\qed \medskip
%%%%%%%%%%%%%%%%%%%%%%%%%%%%%%%%%%%%%%%%%%%%%%%%%%%%%%%%%%%%%%%%%%%%%%

It should be mentioned that, in the case $2\leq n=N$, 
Corollary \ref{c6386} gives bounds for the Euclidean distance 
between the product $\prod_{j=1}^nz_j$  and 
$(\frac{1}{n}\sum_{j=1}^nz_j)^n$.

Let us now discuss the benefit of the bounds 
\eqref{e7465283} and \eqref{e7465284} in the next example. 

%%%%%%%%%%%%%%%%%%%%%%%%%%%%%%%%%%%%%%%%%%%%%%%%%%%%%%%%%%%%%%%%%%%%%%
\begin{example}\label{ex782566}
Let $2\leq n=N$ and let the notation of Theorem \ref{t74739} 
be valid. 
\begin{enumerate}

\item \label{ex782566.a}
We consider the case of $Z\in\{0,1\}^{n\times n}$, where 
$\sum_{\ell=1}^nz_{j,\ell}=\sum_{k=1}^nz_{k,r}=n\widetilde{z}_1$ 
for all $j,r\in\set{n}$, 
i.e.\ the row and column sums of $Z$ are identical.
In particular, $\beta=\widetilde{z}_1^2$. For 
$(u,v)\in\set{n}_{\neq}^2$ and $r\in\set{n}$, we have 
$z_{u,r}^2=z_{u,r}$ and and $y_{u,v,r}^2=|y_{u,v,r}|\in\{0,1\}$. 
Hence 
\begin{align*}
(n(n-1))^3\vartheta^2 
&=\sum_{(r,s)\in \set{n}_{\neq}^2}\Big(\sum_{(u,v)\in\set{n}_{
\neq}^2}y_{u,v,r}^2y_{u,v,s}^2\Big)^2\\
&=\sum_{(r,s)\in \set{n}_{\neq}^2}\Big(\sum_{(u,v)\in\set{n}^2}
(z_{u,r}-2z_{u,r}z_{v,r}+z_{v,r})\,
(z_{u,s}-2z_{u,s}z_{v,s}+z_{v,s})\Big)^2\\
&=\sum_{(r,s)\in \set{n}_{\neq}^2}\Big(\sum_{(u,v)\in\set{n}^2}
(
z_{u,r}z_{u,s}
-z_{u,r}2z_{u,s}z_{v,s}
+z_{u,r}z_{v,s}
-2z_{u,r}z_{v,r}z_{u,s}
\\
&\quad\quad{} 
+2z_{u,r}z_{v,r}2z_{u,s}z_{v,s}
-2z_{u,r}z_{v,r}z_{v,s}
+z_{v,r}z_{u,s}
-z_{v,r}2z_{u,s}z_{v,s}
+z_{v,r}z_{v,s}
)
\Big)^2\\
&=\sum_{(r,s)\in \set{n}_{\neq}^2}\Big(
n(2-8\widetilde{z}_1)\sum_{u\in\set{n}}z_{u,r}z_{u,s}
+4\Big(\sum_{u\in\set{n}}z_{u,r}z_{u,s}\Big)^2
+2n^2\widetilde{z}_1^2\Big)^2,
\end{align*}
from which a formula for $\vartheta$ can be derived. 
Furthermore, for $(u,v),(r,s)\in\set{n}_{\neq}^2$, we have 
\begin{align*}
\lefteqn{\sum_{j\in\set{n}\setminus\{u,v\}}
\sum_{\ell\in\set{n}\setminus\{r,s\}}|z_{j,\ell}|^2
=\sum_{j\in\set{n}\setminus\{u,v\}}
\Big(\sum_{\ell\in\set{n}}z_{j,\ell}-z_{j,r}-z_{j,s}\Big)}\\
&=\sum_{j\in\set{n}}
\Big(\sum_{\ell\in\set{n}}z_{j,\ell}-z_{j,r}-z_{j,s}\Big)
-\Big(\sum_{\ell\in\set{n}}z_{u,\ell}-z_{u,r}-z_{u,s}\Big)
-\Big(\sum_{\ell\in\set{n}}z_{v,\ell}-z_{v,r}-z_{v,s}\Big)\\
&=n^2\widetilde{z}_1-4n\widetilde{z}_1
+z_{u,r}+z_{u,s}+z_{v,r}+z_{v,s}
\leq n^2\widetilde{z}_1-4n\widetilde{z}_1+4
=(n-2)^2\widetilde{z}_1+4(1-\widetilde{z}_1)
\end{align*}
such that 
\begin{align*}
\kappa
\leq \min\Big\{1,
\widetilde{z}_1+\frac{4(1-\widetilde{z}_1)}{(n-2)^2}\Big\}.
\end{align*}
The calculations given here will be used in the subsequent parts of 
this example. 

\item\label{ex782566.b}% 
(Derangement numbers)
Let us now assume that $Z=J-I_n\in\{0,1\}^{n\times n}$, where 
$J$ denotes the matrix all of whose entries are $1$ and $I_n$ 
is the identity matrix. 
Then $\overline{p}_{\set{n}}=\Per(Z)$ 
is the $n$th derangement number, i.e.\ 
the number of permutations in $\set{n}_{\neq}^n$ without fixed 
points. It satisfies
\begin{align*}
\Per(Z)=n!\sum_{j=0}^n\frac{(-1)^j}{j!}
\end{align*}
and can also be interpreted as the number of ways a dance can 
be arranged for $n$ married couples, so that no one dances with 
his or her partner, e.g.\ see Minc \citet[page 44]{MR504978}.
Under the present assumptions, we have 
$\widetilde{z_1}=\frac{n-1}{n}$,
$\beta=\widetilde{z_1}^2=(\frac{n-1}{n})^2$, and
$\sum_{u\in\set{n}}z_{u,r}z_{u,s}=n-2$
for all $(r,s)\in\set{n}_{\neq}^2$. 
Therefore \eqref{ex782566.a} gives 
$\vartheta=\frac{2}{n(n-1)}$ and 
$\kappa\leq\min\{1,\frac{n-1}{n}+\frac{4}{n(n-2)^2}\}$.
Theorems \ref{t74739}, \ref{t4216459} and \eqref{e624373} imply 
that 
\begin{align*} 
\Big|\frac{\overline{p}_{\set{n}}}{n!}-\widetilde{p}_{\set{n}}\Big|
&\leq \frac{f_n(\sqrt{\beta},\sqrt{
\min\{\kappa,\widetilde{\kappa}\}})}{n^2(n-1)}
\leq \frac{1}{2n}.
\end{align*}
In particular, 
$|\frac{\overline{p}_{\set{n}}}{n!}-\widetilde{p}_{\set{n}}|
=O(\frac{1}{n})$ as $n\to\infty$.
Here, using that $(k!)^{1/k}$ is increasing in $k\in\NN$,
it is easily shown that, for $n\geq 4$,
$\widetilde{\kappa}$ from Theorem \ref{t4216459} satisfies 
\begin{align*}
\widetilde{\kappa}
=\frac{1}{n-2}\Big((n-4)\frac{((n-3)!)^{2/(n-3)}}{
((n-2)!)^{2/(n-2)}}+2\Big).
\end{align*}
We note that \eqref{e514385} only gives the bad bound 
$|\frac{\overline{p}_{\set{n}}}{n!}-\widetilde{p}_{\set{n}}|
\leq \frac{n-1}{2n}$. 
Furthermore, the upper bounds in 
\eqref{e4176438}, \eqref{e63286508} and \eqref{e4176439} from the 
introduction cannot be small, since 
$\gamma=\frac{n-1}{2n-1}$, where $\gamma$ is defined there. 

\item \label{ex782566.c}
(M\'{e}nage numbers)
We now consider the matrix $Z=J-I_n-P\in\{0,1\}^{n\times n}$, 
where $n\geq 3$ and $P$ is the matrix with $1$'s in positions 
$(1,2),\,(2,3),\dots,(n-1,n),\,(n,1)$ and $0$'s otherwise. 
Then $\overline{p}_{\set{n}}=\Per(Z)$ 
is the $n$th m\'{e}nage number, which can be described as the
number of ways of seating a set of 
married couples at a circular table so that men and women 
alternate and nobody sits next to his or her partner, 
see Minc \citet[page 44]{MR504978}. 
The following explicit formula is due to Touchard 
\citet{Touchard1934}: 
\begin{align*} 
\Per(Z)= \sum_{j=0}^n(-1)^j\frac{2n}{2n-j}\binomial{2n-j}{j}(n-j)!.
\end{align*}
In this situation, we have  $\widetilde{z_1}=\frac{n-2}{n}$,
$\beta=\widetilde{z_1}^2=(\frac{n-2}{n})^2$,
and, for $(r,s)\in\set{n}_{\neq}^2$, 
\begin{align*}
\sum_{u\in\set{n}}z_{u,r}z_{u,s}
=\begin{cases}
n-3,  & \mbox{ if } 
\card{\{(r,s),(s,r)\}\cap\{(1,2),\dots,(n-1,n),(n,1)\}}=1,\\
n-4,  & \mbox{ otherwise}.
\end{cases}
\end{align*}
Therefore \eqref{ex782566.a} gives 
$\vartheta=\frac{\sqrt{8(n^2+4n-20)}}{n(n-1)^{3/2}}$ and 
$\kappa\leq\min\{1,\frac{n-2}{n}+\frac{8}{n(n-2)^2}\}$.
Theorems \ref{t74739}, \ref{t4216459} and \eqref{e624373} imply
that
\begin{align*}  
\Big|\frac{\overline{p}_{\set{n}}}{n!}
-\widetilde{p}_{\set{n}}\Big|
&\leq\frac{\sqrt{2(n^2+4n-20)}}{n^{2}(n-1)^{3/2}}
f_n(\sqrt{\beta},\sqrt{\min\{\kappa,\widetilde{\kappa}\}})
\leq\frac{\sqrt{n^2+4n-20}}{\sqrt{2(n-1)\,}n},
\end{align*}
i.e.\ 
$|\frac{\overline{p}_{\set{n}}}{n!}-\widetilde{p}_{\set{n}}|
=O(\frac{1}{\sqrt{n}})$, as $n\to\infty$. 
As in Part \eqref{ex782566.b}, \eqref{e514385} only gives a bad 
bound, namely
$|\frac{\overline{p}_{\set{n}}}{n!}-\widetilde{p}_{\set{n}}|
\leq \frac{n-2}{n}$.
Again, the bounds in \eqref{e4176438}, \eqref{e63286508} 
and \eqref{e4176439} from the introduction cannot be small, since 
$\gamma=\frac{n-2}{2(n-1)}$.
We finally note that, for $n\geq 5$, $\widetilde{\kappa}$ can 
easily be evaluated as
\begin{align*}
\widetilde{\kappa}=\frac{1}{n-2}
\Big((n-5)\frac{((n-4)!)^{2/(n-4)}}{((n-2)!)^{2/(n-2)}}
+2\frac{((n-3)!)^{2/(n-3)}}{((n-2)!)^{2/(n-2)}}+1\Big).
\end{align*}
\end{enumerate}
\end{example}
%%%%%%%%%%%%%%%%%%%%%%%%%%%%%%%%%%%%%%%%%%%%%%%%%%%%%%%%%%%%%%%%%%%%%%

\subsection{Second order approximation} 
%%%%%%%%%%%%%%%%%%%%%%%%%%%%%%%%%%%%%%%%%%%%%%%%%%%%%%%%%%%%%%%%%%%%%%

The next theorem contains an improvement of inequality 
\eqref{e4176439} from the introduction. 
%%%%%%%%%%%%%%%%%%%%%%%%%%%%%%%%%%%%%%%%%%%%%%%%%%%%%%%%%%%%%%%%%%%%%%
\begin{theorem}\label{th487698}
Let $2\leq n\leq N$,
$\widetilde{p}^{(2)}
=\sum_{R\subseteq\set{n}:\,|R|=2} 
\widetilde{p}_{\set{n}\setminus R}
\sum_{j=1}^N\prod_{r\in R}(z_{j,r}-\widetilde{z}_r)$
and $\beta=\frac{1}{n}\sum_{r=1}^n|\widetilde{z}_r|^2$. 
Let 
\begin{gather*}
\vartheta_3
=\frac{(N-3)!}{N!}\sqrt{\frac{(n-3)!}{n!}}
\Big(\sum_{(r,s,t)\in \set{n}_{\neq}^3}
\Big(\sum_{(u,v,w)\in\set{N}_{\neq}^3}|y_{u,v,r}y_{u,v,s}y_{u,w,t}|
\Big)^2\Big)^{1/2},
\quad\mbox{ if } n\geq 3,\\
\vartheta_4
=\frac{(N-4)!}{N!}\sqrt{\frac{(n-4)!}{n!}}
\Big(\sum_{(q,r,s,t)\in \set{n}_{\neq}^4}
\Big(\sum_{(u,v,w,x)\in\set{N}_{\neq}^4}
|y_{u,v,q}y_{u,v,r}y_{w,x,s}y_{w,x,t}|
\Big)^2\Big)^{1/2},\quad\mbox{ if } n\geq 4.
\end{gather*}
If $n=2$, then $\vartheta_3=0$; 
further, if $n\in\{2,3\}$, then $\vartheta_4=0$.
Let 
\begin{gather*}
\kappa^{(\nu)}
= 
\begin{cases}
\displaystyle
\frac{1}{(n-\nu)(N-\nu)}
\max_{J\subseteq\set{N}, \, R\subseteq\set{n}:\,
\card{J}=\card{R}=\nu}
\sum_{j\in\set{N}\setminus J}
\sum_{r\in\set{n}\setminus R}|z_{j,r}|^2, &
\mbox{if } n\geq \nu+1,\\
1, &\mbox{if } n\leq \nu,
\end{cases}\quad(\nu\in\{3,4\})\\
f_n(x_1,x_2)
=\sum_{k\in\set{n}\setminus\set{2}}
(n+k-2)(n-k+1)x_1^{n-k}x_2^{k-3},\\
g_n(x_1,x_2)
=\sum_{k\in\set{n}\setminus\set{3}}
(k-3)(n+k-2)(n-k+1)x_1^{n-k}x_2^{k-4},
\end{gather*}
for $x_1,x_2\in\RR$. Then
\begin{align}
\Bigl|\frac{(N-n)!}{N!}\overline{p}_{\set{n}}
-\widetilde{p}_{\set{n}}
+\frac{\widetilde{p}^{(2)}}{N(N-1)}\Bigr|
&\leq\frac{\vartheta_3}{2N^2} f_n(\sqrt{\beta},\sqrt{\kappa^{(3)}})
+\frac{\vartheta_4}{8N^2} g_n(\sqrt{\beta},\sqrt{\kappa^{(4)}}).
\label{e5196573}
\end{align}
\end{theorem}
%%%%%%%%%%%%%%%%%%%%%%%%%%%%%%%%%%%%%%%%%%%%%%%%%%%%%%%%%%%%%%%%%%%%%%
\Proof 
Theorem \ref{t384696}
implies that 
\begin{align*}
\Bigl|\frac{(N-n)!}{N!}\overline{p}_{\set{n}}
-\widetilde{p}_{\set{n}}
+\frac{\widetilde{p}^{(2)}}{N(N-1)}\Bigr|
\leq \frac{M_1}{2N^2}+\frac{M_2}{8N^2},
\end{align*}
where
\begin{align*}
M_1&=\frac{(N-n)!}{N!}\sum_{k\in\set{n}\setminus\set{2}}h_{k,n}
\sum_{\newatop{R\subseteq\set{n}}{\card{R}=k}}
\sum_{(r,s,t)\in R_{\neq}^3}
|\widetilde{p}_{\set{n}\setminus R}|
\sum_{j\in\set{N}_{\neq}^n}|y_{j_r,j_s,r}y_{j_r,j_s,s}y_{j_r,j_t,t}
p_{j,R\setminus\{r,s,t\}}|,\\
M_2&=\frac{(N-n)!}{N!}\sum_{k\in\set{n}\setminus\set{3}}h_{k,n}
\sum_{\newatop{R\subseteq \set{n}}{\card{R}=k}}
\sum_{(q,r,s,t)\in R_{\neq}^4}
|\widetilde{p}_{\set{n}\setminus R}|
\sum_{j\in\set{N}_{\neq}^n}
|y_{j_q,j_r,q}y_{j_q,j_r,r}y_{j_s,j_t,s}y_{j_s,j_t,t}
p_{j,R\setminus\{q,r,s,t\}}|.
\end{align*}
Similarly as in the proof of Theorem \ref{t74739},
one can apply Lemma \ref{l378687}\eqref{l378687.a},
the Cauchy-Schwarz inequality in combination with 
Maclaurin's inequality to show that 
$M_1\leq \vartheta_3f_n(\sqrt{\beta},\sqrt{\kappa^{(3)}})$ and 
$M_2\leq\vartheta_4g_n(\sqrt{\beta},\sqrt{\kappa^{(4)}})$.
\hfill\qed \medskip

%%%%%%%%%%%%%%%%%%%%%%%%%%%%%%%%%%%%%%%%%%%%%%%%%%%%%%%%%%%%%%%%%%%%%%

A theorem for $(0,1)$-matrices $Z$ similar to the above can be 
proved with the help of Lemma~\ref{l378687}\eqref{l378687.b}. 
Further, 
the auxiliary inequalities contained in the next lemma can be used in 
combination with \eqref{e624372} to prove upper bounds of the 
right-hand side of \eqref{e5196573}. For instance, it is possible 
to give an estimate, which is of the same order as the right-hand 
side of \eqref{e4176439}, if $\gamma$ is bounded away from~$1$. 
Since \eqref{e624372} and \eqref{e8276436} are based on the 
Cauchy-Schwarz inequality, the form of \eqref{e5196573} is better 
than that of \eqref{e4176439}. We omit the details here. 
%%%%%%%%%%%%%%%%%%%%%%%%%%%%%%%%%%%%%%%%%%%%%%%%%%%%%%%%%%%%%%%%%%%%%%

\begin{lemma} 
Let us assume that $2\leq n\leq N$. Let 
$\beta$, $f_n$ and $g_n$ be as in Theorem \ref{th487698}.
Let $\vartheta$ as in Theorem \ref{t74739}. Then, for $x\in[0,1]$, 
\begin{gather}
\vartheta_3 
\leq \frac{(N-2)!}{N!}\sqrt{\frac{(n-3)!}{n!}}
\sum_{(u,v)\in\set{N}_{\neq}^3}
\Big(\sum_{r\in \set{n}}|y_{u,v,r}|^2\Big)^{3/2},\label{e8276436}\\
\vartheta_4 
\leq \sqrt{\frac{n(n-1)}{(n-2)(n-3)}}
\frac{N(N-1)}{(N-2)(N-3)}\vartheta^2,\label{e8276437}\\
f_n(x,1)
\leq 2(n-1)\min\Big\{\frac{n(n-2)}{3},\frac{1}{(1-x)^2}\Big\}
,\label{e765143255}\\
g_n(x,1)
\leq2(n-1)(n-3)\min\Big\{\frac{n(n-2)}{8},\frac{1}{(1-x)^2}\Big\}.
\label{e768457}
\end{gather}
\end{lemma}
%%%%%%%%%%%%%%%%%%%%%%%%%%%%%%%%%%%%%%%%%%%%%%%%%%%%%%%%%%%%%%%%%%%%%%
\Proof
Using that 
\begin{align*}
\lefteqn{\sum_{(r,s,t)\in \set{n}_{\neq}^3}
\Big(\sum_{(u,v,w)\in\set{N}_{\neq}^3}|y_{u,v,r}y_{u,v,s}y_{u,w,t}|
\Big)^2}\\
&\qquad=\sum_{(u,v,w)\in\set{N}_{\neq}^3}
\sum_{(u',v',w')\in\set{N}_{\neq}^3}
\sum_{(r,s,t)\in \set{n}_{\neq}^3}
|y_{u,v,r}y_{u,v,s}y_{u,w,t}
y_{u',v',r}y_{u',v',s}y_{u',w',t}|
\end{align*}
together with the Cauchy-Schwarz inequality, we get 
\begin{align*}
\Big(\sum_{(r,s,t)\in \set{n}_{\neq}^3}
\Big(\sum_{(u,v,w)\in\set{N}_{\neq}^3}|y_{u,v,r}y_{u,v,s}y_{u,w,t}|
\Big)^2\Big)^{1/2}
\leq \sum_{(u,v,w)\in\set{N}_{\neq}^3}
\Big(\sum_{(r,s,t)\in \set{n}_{\neq}^3}
|y_{u,v,r}y_{u,v,s}y_{u,w,t}|^2\Big)^{1/2}
\end{align*}
and this, in turn, can be estimated with the help of the H\"older 
inequality by 
\begin{align*}
\lefteqn{\sum_{(u,v,w)\in\set{N}_{\neq}^3}
\Big(\sum_{r\in \set{n}}|y_{u,v,r}|^2\Big)
\Big(\sum_{r\in \set{n}}|y_{u,w,r}|^2\Big)^{1/2}}\\
&\leq \Big(\sum_{(u,v,w)\in\set{N}_{\neq}^3}
\Big(\sum_{r\in \set{n}}|y_{u,v,r}|^2\Big)^{3/2}\Big)^{2/3}
\Big(\sum_{(u,v,w)\in\set{N}_{\neq}^3}
\Big(\sum_{r\in \set{n}}|y_{u,w,r}|^2\Big)^{3/2}\Big)^{1/3}\\
&=(N-2)\sum_{(u,v)\in\set{N}_{\neq}^3}
\Big(\sum_{r\in \set{n}}|y_{u,v,r}|^2\Big)^{3/2},
\end{align*}
which implies \eqref{e8276436}. 
Further, 
\begin{align*}
\lefteqn{
\sum_{(q,r,s,t)\in \set{n}_{\neq}^4}
\Big(\sum_{(u,v,w,x)\in\set{N}_{\neq}^4}
|y_{u,v,q}y_{u,v,r}y_{w,x,s}y_{w,x,t}|
\Big)^2}\\
&\leq \sum_{(q,r)\in \set{n}_{\neq}^2}
\sum_{(s,t)\in \set{n}_{\neq}^2}
\Big(\sum_{(u,v)\in\set{N}_{\neq}^2}
|y_{u,v,q}y_{u,v,r}|
\sum_{(w,x)\in\set{N}_{\neq}^2}
|y_{w,x,s}y_{w,x,t}|\Big)^2\\
&=\Big(\sum_{(q,r)\in \set{n}_{\neq}^2}
\Big(\sum_{(u,v)\in\set{N}_{\neq}^2}
|y_{u,v,q}y_{u,v,r}|\Big)^2\Big)^2,
\end{align*}
from which  \eqref{e8276437} follows. 
For \eqref{e765143255}, we note that
\begin{align*}
f_n(x,1)
&\leq \sum_{k\in\set{n}\setminus\set{2}}(n+k-2)(n-k+1) 
=\frac{2}{3}n(n-1)(n-2)
\end{align*}
and $f_n(x,1)
=\sum_{k\in\set{n}\setminus\set{2}}(n+k-2)(n-k+1)x^{n-k}
\leq 2(n-1)\sum_{k\in\set{n}\setminus\set{2}}(n-k+1)x^{n-k}$,
where
$\sum_{k\in\set{n}\setminus\set{2}}(n-k+1)x^{n-k}
\leq \sum_{k=0}^{\infty}(k+1)x^{k}
=\frac{1}{(1-x)^2}$.
Furthermore, 
\begin{align*}
g_n(x,1)
=\sum_{k\in\set{n}\setminus\set{3}}(k-3)(n+k-2)(n-k+1)x^{n-k}
\leq (n-3)f_n(x,1)
\leq \frac{2(n-1)(n-3)}{(1-x)^2}
\end{align*}
and 
$g_n(x,1)
\leq\sum_{k\in\set{n}\setminus\set{2}}(k-3)(n+k-2)(n-k+1)
=\frac{n!}{4(n-4)!}$, which implies \eqref{e768457}.\hfill \qed
%%%%%%%%%%%%%%%%%%%%%%%%%%%%%%%%%%%%%%%%%%%%%%%%%%%%%%%%%%%%%%%%%%%%%%
\section*{Acknowledgment}
%%%%%%%%%%%%%%%%%%%%%%%%%%%%%%%%%%%%%%%%%%%%%%%%%%%%%%%%%%%%%%%%%%%%%%
The author thanks Lingji Chen for helpful discussions and for pointing
out a flaw in a previous version of Theorem \ref{th487698}.
%%%%%%%%%%%%%%%%%%%%%%%%%%%%%%%%%%%%%%%%%%%%%%%%%%%%%%%%%%%%%%%%%%%%%%

\small 
\raggedbottom % inserted to avoid Underfull \vbox (badness 10000)
\linespread{1}
\selectfont
\bibliography{npa97_ber}

%%%%%%%%%%%%%%%%%%%%%%%%%%%%%%%%%%%%%%%%%%%%%%%%%%%%%%%%%%%%%%%%%%%%%%
\end{document}